\documentclass[12pt]{article} 
\usepackage{hyperref}
\usepackage{array,epsfig,fancyheadings,rotating}

\usepackage{sectsty, secdot}
\sectionfont{\fontsize{12}{14pt plus.8pt minus .6pt}\selectfont}
\renewcommand{\theequation}{\thesection\arabic{equation}}
\subsectionfont{\fontsize{12}{14pt plus.8pt minus .6pt}\selectfont}

\usepackage[T1]{fontenc}        

\usepackage[utf8x]{inputenc}    
\usepackage[comma,authoryear]{natbib}
\usepackage{ucs}                
\usepackage{amsmath}            
\usepackage{amsfonts}
\usepackage{amssymb}
\usepackage{amsthm} 
\usepackage{thmtools}
\usepackage{ifthen}
\usepackage{color}
\usepackage{graphicx} 
\usepackage{textcomp}            
\usepackage{palatino}            
\usepackage{eulervm}            
\usepackage{varioref}            
\usepackage{cleveref}             
\usepackage{float}
\usepackage{multirow}
\usepackage{rotating}
\usepackage{subfig}
\usepackage{xr} 
\usepackage{setspace} 

\newcommand{\R}{{\mathbb{R}}}         
\newcommand{\E}{{\mathbb{E}}}
 \newcommand{\N}{{\mathbb{N}}}         
\newcommand{\bqan}{\begin{eqnarray}}
\newcommand{\eqan}{\end{eqnarray}}
\newcommand{\bit}{\begin{itemize}}
\newcommand{\eit}{\end{itemize}}

\newcommand{\vepsilon}{\vep}

\newcommand{\Cov}{\operatorname{{\textit Cov}}}
\newcommand{\var}{\operatorname{{\textit Var}}}
\newcommand{\Var}{\operatorname{{\textit Var}}}

\newcommand{\diag}{\operatorname{diag}}

\newcommand{\rank}{\operatorname{\textit rank}}

\newcommand{\vh}{\boldsymbol{h}}
\newcommand{\vn}{\boldsymbol{n}}
\newcommand{\vv}{\boldsymbol{v}}

\newcommand{\vC}{\boldsymbol{C}}

\newcommand{\vE}{\boldsymbol{E}}
\newcommand{\vI}{\boldsymbol{I}}
\newcommand{\vJ}{\boldsymbol{J}}
\newcommand{\vP}{\boldsymbol{P}}
\newcommand{\vV}{\boldsymbol{V}}
\newcommand{\vX}{\boldsymbol{X}}
\newcommand{\vY}{\boldsymbol{Y}}
\newcommand{\vZ}{\boldsymbol{Z}}
\newcommand{\vN}{\boldsymbol{N}}
\newcommand{\vDelta}{\boldsymbol{\Delta}}
\newcommand{\vep}{\boldsymbol{\epsilon}}
\newcommand{\valpha}{\boldsymbol{\alpha}}
\newcommand{\vbeta}{\boldsymbol{\beta}}
\newcommand{\vgamma}{\boldsymbol{\gamma}}
\newcommand{\vmu}{\boldsymbol{\mu}}
\newcommand{\vnu}{\boldsymbol{\nu}}
\newcommand{\vSigma}{\boldsymbol{\Sigma}}
\newcommand{\vzeta}{\boldsymbol{\zeta}}
\newcommand{\veins}{{\bf 1}}
\newcommand{\vnull}{{\bf 0}}
\newcommand{\ve}{\boldsymbol{e}}
\newcommand{\To}{\longrightarrow}            
\newcommand{\ind}{1\hspace{-0.7ex}1}

\DeclareMathOperator{\vech}{vech}
\DeclareMathOperator{\tr}{tr}

\newtheorem{theorem}{Theorem}

\newtheorem{corollary}{Corollary}

\theoremstyle{definition}

\newtheorem{remark}{Remark}
\newcommand{\lan}{ \scriptstyle \mathcal{O}\textstyle}

\makeatletter
\renewenvironment{proof}[1][\proofname]{\par
  \pushQED{\qed}%
  \fontfamily{ppl}\fontseries{m}\fontshape{it} \topsep6\p@\@plus6\p@\relax
  \trivlist
  \item[\hskip\labelsep
        \bfseries
    #1\@addpunct{:}]\ignorespaces
}{%
  \popQED\endtrivlist\@endpefalse
}
\makeatother

\relax 
\begin{document}
\renewcommand{\baselinestretch}{2}

\markright{ \hbox{\footnotesize\rm Statistica Sinica
}\hfill\\[-13pt]
\hbox{\footnotesize\rm
}\hfill }

\markboth{\hfill{\footnotesize\rm Paavo Sattler, Arne Bathke and Markus Pauly} \hfill}
{\hfill {\footnotesize\rm FILL IN A SHORT RUNNING TITLE} \hfill}

\renewcommand{\thefootnote}{}
$\ $\par


\fontsize{12}{14pt plus.8pt minus .6pt}\selectfont \vspace{0.8pc}
\centerline{\large\bf Testing  Hypotheses about }
\vspace{2pt} 
\centerline{\large\bf    Covariance Matrices in General MANOVA Designs}
\vspace{.4cm} 
\centerline{Paavo Sattler$^1$, Arne Bathke$^2$ and Markus Pauly$^1$} 
\vspace{.4cm} 
\centerline{\it ${}^1$TU Dortmund University and ${}^2$University of Salzburg}
 \vspace{.55cm} \fontsize{9}{11.5pt plus.8pt minus.6pt}\selectfont

\begin{quotation}
\noindent 
\begin{abstract}
\noindent We introduce a unified approach to testing a variety of rather general null hypotheses that can be formulated in terms of covariances matrices. These include as special cases, for example, testing for equal variances, equal traces, or for elements of the covariance matrix taking certain values. 
The proposed method only requires very few assumptions and thus promises to be of broad practical use. 
Two test statistics are defined, and their asymptotic or approximate sampling distributions are derived.
In order to particularly improve the small-sample behavior of the resulting tests, two bootstrap-based methods are developed and theoretically justified. 
Several simulations shed light on the performance of the proposed tests. 
The analysis of a real data set illustrates the application of the procedures.
\end{abstract}


\vfill
\vfill

\newpage
\par
\end{quotation}\par

\def\thefigure{\arabic{figure}}
\def\thetable{\arabic{table}}

\renewcommand{\theequation}{\thesection.\arabic{equation}}

\fontsize{12}{14pt plus.8pt minus .6pt}\selectfont
\section{MOTIVATION AND INTRODUCTION}\label{int}

It is of substantial interest to have valid statistical methods for inference on covariance matrices available for at least two major reasons.
The first one is that a treatment effect may indeed best be described by a particular configuration of scale or covariance parameters -- not by a mean difference. The second reason corresponds to a more indirect purpose, namely that the main interest of the investigation may be described by a location change under the alternative, but some of the available inference methods for location effects rely on assumptions regarding variances or covariances that need to be assessed reliably. In either situation, a statistical test about hypotheses that are formulated in terms of covariance matrices is necessary. From a methodological point of view, such a test shall not make too many restrictive assumptions itself, for example, regarding underlying distributions. Furthermore, it shall perform well for moderate sample sizes, where clearly the term {\em moderate} will have to be seen in connection with the number of parameters effectively being tested.

Considering the central importance and the widespread need for hypothesis tests on covariance matrices, it may come as a surprise that a general and unifying approach to this task has not been developed thus far. There are several tests for specialized situations, such as testing equality of variances or even covariance matrices. Many of these approaches will be mentioned below. 
However, they typically only address one particular question, and they often rely on restrictive distributional assumptions,
such as normality (e.g.~in \cite{box1953} and \cite{anderson1984}), elliptical distributions (e.g. in \cite{muirhead}, \cite{fangZhang} and \cite{hallin}), or conditions on the characteristic functions (e.g. in \cite{gupta2006}).

One exception is the test of \cite{zhang1993}, which theoretically allows for testing a multitude of hypotheses without restrictive distributional conditions. Unfortunately, this procedure's small and medium sample performance is comparatively poor, particularly regarding the power. Their technique to improve the performance requires a more restrictive null hypothesis that additionally postulates equality of certain moments. This makes it somewhat difficult to use this approach in practice, as rejection does not mean that the covariances are unequal.

The goal of the present article is to introduce a very general approach to statistical hypothesis testing, where the hypotheses are formulated in terms of covariance matrices. This includes as special cases, for example, hypotheses formulated using their traces, hypotheses of equality of variances or of covariance matrices, and hypotheses in which a covariance matrix is assumed to have particular entries. 
The test procedures are based on a resampling approach whose asymptotic validity is shown theoretically, while the actual finite sample performance has been investigated through extensive simulation studies. Analysis of a real data example illustrates the application of the proposed methods. 

In the following section, the statistical model and (examples for) different null hypotheses that can be investigated using the proposed approach will be introduced. Thereafter, the asymptotic distributions of the proposed test statistics are derived (Section~\ref{The Test Statistics and its Asymptotics}) and proven to be regained by two different resampling strategies (Section~\ref{Resampling Procedures}). The simulation results regarding type-I-error control and power are discussed in Section~\ref{Simulations}, computation time is considered in Section~\ref{Time}, while an illustrative data analysis of EEG-data is conducted in Section~\ref{Illustrative Data Analysis}. All proofs are deferred to a technical supplement.

\section{\textsc{Statistical Model and Hypotheses}} \label{mod}

We consider a general semiparametric model given by independent $d$-dimensional random vectors \bqan\label{model}
\vX_{ik}&= \vmu_i + \vep_{ik}. 
\eqan
Here, the index $i=1, \dots, a$ refers to the treatment group and $k=1, \dots,n_i$ to the individual, on which $d$-dimensional observations are measured. More details to this model can be found in the supplementary material.\\
In this setting, $\E(\vX_{ik})=\vmu_i = (\mu_{i1}, \dots, \mu_{id})^\top \in \mathbb{R}^d$ denotes the $i$-th group mean while 
the residuals $\vep_{i1},\dots,\vep_{in_i}$ are assumed to be centered $\E(\vep_{i1}) = \vnull$ and i.i.d. within each group. We require 
finite fourth moment $\E(||\vep_{i1}||^4) < \infty$, where this denotes the Euclidean norm. However, beyond this, there are no other distributional assumptions. In particular, the covariance matrices $\Cov(\vep_{i1})=\vV_i \geq 0$ may be arbitrary and do not even have to be positive definite. 
For convenience, we aggregate the individual vectors into $\vX=(\vX_{11}^\top, \dots, \vX_{an_a}^\top)^\top$ as well as $\vmu = (\vmu_{1}^\top, \dots, \vmu_a^\top)^\top$. Stacking the covariance matrix $\vV_i =(v_{irs})_{r,s}^d$ into the $p:=d(d+1)/2$-dimensional vector $\vv_i = \vech(\vV_i) = (v_{i11},v_{i12},\dots,v_{i1d},v_{i22},\dots,v_{i2d},\dots,v_{idd})^\top$  $(i=1,\dots,a)$ containing  the upper triangular entries of $\vV_i$  we formulate hypotheses in terms of the pooled covariance vector 
$\vv = (\vv_1^\top,\dots,\vv_a^\top)^\top$ as 
\bqan\label{eq:hypo cov}
 \mathcal H_0^{\vv}: \vC \vv = \vzeta.
\eqan
Here, $\vC$ denotes a suitable hypothesis matrix of interest, and $\vzeta$ is a fixed vector. It should be noted that we don't assume that $\vC$ is a contrast matrix, not to mention a projection matrix. This is different from the frequently used hypothesis formulation about mean vectors in MANOVA designs (\cite{kon:2015}, \citet*{friedrich2017permuting}, \cite{bathke2018}), where one can usually work with a unique projection matrix. 
However, working with simpler matrices (as we do)  can help to save considerable computation time, see Remark~\ref{Remark Comp Time} below.

In order to discuss some particular hypotheses included within the general setup \eqref{eq:hypo cov}, we fix the following notation: Let $\vI_d$ be the $d$-dimensional unit matrix, $\boldsymbol{1}_d=(1, \dots, 1)^\top$ the $d$-dimensional column vector of 1's and $\vJ_d = \boldsymbol{1}_d \boldsymbol{1}_d^\top$ the $d$-dimensional matrix of 1’s. Furthermore, $\vP_a = \vI_a -  \vJ_a/a$ denotes the $a$-dimensional centering matrix, while $\oplus$ and $\otimes$ denote direct sum and Kronecker product, respectively.
Then the following null hypotheses of interest are covered:

(a) {\bf Testing equality of variances:} For a univariate outcome with $d=1$, testing the null hypothesis 
 $
 \mathcal H_0^{\vv}: v_{1 11} = v_{211} = \dots = v_{a11}
 $
 of equal variances is included within \eqref{eq:hypo cov} by setting 
 $\vC = \vP_a$ and $\vzeta=\vnull$. Hypotheses of this type have been studied by \cite{bartlett1953} as well as \cite{boos2004}, \cite{gupta2006}, and \cite{pauly:2011a}, among others. In the special case of a two-armed design with $a=2$,  this is also the null hypothesis inferred by the popular F-ratio test which, however, is known to be sensitive to deviations from normality (\cite{box1953}).
 
 (b)  {\bf Testing for a given covariance matrix:} Let $\vSigma$ be a given covariance matrix. It may represent, for example, an autoregressive or compound symmetry covariance structure. For $a=1$, our general formulation also covers testing the null hypothesis 
 $
 \mathcal  H_0^{\vv}: \vV_1 = \vSigma
 $
 by setting $\vC = \vI_p$ and $\vzeta = \vech(\vSigma)$. 
 Hypotheses of this kind have been studied by \cite{gupta2006}.
 
 (c)  {\bf Testing homogeneity of covariance matrices:} More general than in (a), let $\vC = \vP_a\otimes \vI_p$ and $\vzeta=\vnull$ for arbitrary $d\in\N$. Then \eqref{eq:hypo cov} describes the null hypothesis
 $ H_0^{\vv}: \vV_1 = \dots = \vV_a.
 $
 For multivariate normally distributed random variables, this is the testing problem of Box's-M-test \cite{box1953}, for which extensions have been studied in \cite{lawley1963}, \cite{browne}, \citet*{zhu}, and \cite{yang}. Moreover, \cite{zhang1992, zhang1993} proposed Bartlett-type tests with bootstrap approximations in a general model similar to ours. However, the pooled bootstrap method of \cite{zhang1992} requires equality of some special kind of fourth moments across groups while the separate bootstrap approximation proposed in \cite{zhang1993} exhibited unsatisfactory small sample behavior in terms of size control or power. \\

 Beyond the above choices, $\mathcal H_0^{\vv}$ in \eqref{eq:hypo cov} even contains hypotheses about linear functions of matrices. 
 To this end, set $\vh_d:=(1,\vnull_{d-1}^\top, 1, \vnull_{d-2}^\top, \dots, 1,0,1)^\top$ and consider the following examples:
 
 (d)  {\bf Traces as effect measures:} Suppose we are interested in the total variance $\sum_{\ell=1}^d \var(X_{i1\ell}) = \tr (\vV_i)$ of all components as a univariate effect measure for each group. 
This may be an advantageous approach in terms of power, as illustrated in the data example analysis below. Then, their equality  $
  \mathcal H_0^{\vv}:\tr( \vV_1) = \dots = \tr(\vV_a)
 $
 can be tested by choosing $\vC = \vP_a\otimes [\vh_d\cdot \vh_d^\top]/d, $ and $\vzeta=\vnull$

 (e)  {\bf Testing for a given trace:} Consider the situation of example (d) with just one group $a=1$. We then may be interested in testing for a given value $\gamma\in \R$ of the trace, i.e.
  $
  \mathcal H_0^{\vv}: \tr(\vV_1) = \gamma.
 $ {Therefore we chose $\vC=\ve_1\cdot\vh_d^\top$ and $\vzeta=\ve_1\cdot \gamma$, with $\ve_1=(1,\vnull_{d-1}^\top)^\top$.}
 
(f)  {\bf Higher Way Layouts:} Moreover, we can even infer hypotheses about variances, covariance matrices, or traces in arbitrarily crossed multivariate layouts by splitting up indices. For example, consider a two-way cross-classified design with fixed factors $A$ and $B$ whose levels are $i_1 = 1,\dots, a$ and $i_2=1,\dots,b$, respectively. 
Assume that the interest lies in measuring, for example, their effect on the total variance, that is, the trace (a similar approach works for variances and covariances). 
We observe $n_{i_1i_2} > 0$ subjects for each factor level combination 
 $(i_1, i_2)$.  To formulate hypotheses of no main trace effects for each factor, as well as hypotheses of no interaction trace effects we write $\tr(\vV_{i_1i_2}) = t + \alpha_{i_1} + \beta_{i_2} + (\alpha\beta)_{i_1i_2}$ with the usual side conditions $\sum_{i_1}\alpha_{i_1} = \sum_{i_2}\beta_{i_2} = \sum_{i_1}(\alpha\beta)_{i_1} = \sum_{i_2}(\alpha\beta)_{i_2} = 0$. 
Here, for example, $\alpha_{i_1}$ can be interpreted as the part of the total variance under factor level $i_1$ by factor A.  Then, the choice  $\vC=(\vP_a\otimes \vJ_b/b)\otimes(\vh_d\cdot\vh_d^\top/d)$ and $\vzeta=\vnull$ leads to a test for no main effect of factor $A$ (measured in the above trace effects),
 $
 \mathcal H_0^{\vv}: \alpha_{1} = \dots = \alpha_a = 0,
 $
  while $\vC=(\vP_a\otimes \vP_b)\otimes(\vh_d\cdot\vh_d^\top/d)$  and $\vzeta=\vnull$ result in the hypothesis of no interaction (again measured in trace effects) between the factors $A$ and $B$,
  $
  \mathcal H_0^{\vv}: \alpha\beta_{ij} \equiv 0 \text{ for all } i,j.
  $\begin{remark}\label{Remark Comp Time}
Although in most of the considered scenarios, it is possible to find an idempotent symmetric hypothesis matrix $\vC$, 
the option {$\vzeta\neq \vnull_p$} allows for matrices that are neither symmetric nor idempotent. 
From a theoretical point of view, this does not really matter. 
However, from a practical point of view, the choice of the hypothesis matrix may actually have a great effect with regard to saving computation time. To this aim, we allow {$\vC\in\R^{m \times ap}$ with $m\leq ap$ together with appropriate $\vzeta\in \R^m$ and formulate all our theorems for this kind of matrices}. For  example $\mathcal H_0^{\vv}:\tr( \vV_1) =\gamma$ could also be formulated by $\vh_d^\top\cdot \vv=\gamma$. 
{Depending on the hypothesis of interest, the computational savings in our simulations were up to 66\% for smaller dimensions and partially even more than 99\% for larger dimensions, see \Cref{Hmatrix} for a detailed discussion.}  

\end{remark}
 In the subsequent sections, we develop testing procedures for $\mathcal H_0^{\vv}$ in \eqref{eq:hypo cov} and thus for all given examples  (a)--(f) above. The basic idea is to use a quadratic form in the vector $\vC \widehat{\vv} - \vzeta$ of estimated and centered effects. For ease of presentation and its widespread use in our setting (with $\E(||\vep_{i1}||^4) < \infty$), we thereby focus on empirical covariance matrices 
 $
 \widehat{\vV}_i = ({n_i-1})^{-1} \sum_{k=1}^{n_i} (\vX_{ik} - \overline \vX_{i\cdot}) (\vX_{ik} - \overline \vX_{i\cdot})^\top, \quad \widehat{\vv}_i=\vech(\widehat \vV_i), 
 $
 as estimators for $\vV_i$, $i=1,...,a,$ where 
 $\overline \vX_{i\cdot}={n_i}^{-1} \sum_{k=1}^{n_i} \vX_{i k}$.  
Other choices, as, for example, surveyed in \citet*{Duembgen}, may be part of future research.
 
Thereby, inverting the resulting test procedures will lead to \emph{confidence regions} about the effect measures of interest. For example, in case (e), we may obtain confidence intervals for the unknown trace $\tr(\vV_1)$. 

\section{\textsc{The Test Statistics and their Asymptotics}}\label{The Test Statistics and its Asymptotics}

In order to obtain the mentioned inference procedures which are formulated using quadratic forms, we first have to study 
the asymptotic behaviour of the normalized $m$-dimensional vector $\sqrt{N} (\vC\widehat \vv -\vzeta)$, where 
$\widehat \vv = (\widehat{\vv}_1^\top,\dots, \widehat{\vv}_a^\top)^\top$ is the pooled empirical covariance estimator of $\vv$.
For convenience, we thereby assume throughout that the following asymptotic sample size condition holds, as $\min(n_1,\dots,n_a)\to\infty$:

\begin{center}
\bit
  \item[(A1)] $\frac{n_i}{N}\to \kappa_i\in (0,1],~i= 1,...,a$ for $N = \sum_{i=1}^a n_i$. 
\eit
\end{center}
As  $\kappa_i >0$ holds for all $i$, we have $\kappa_1=1$ if and only if $a=1$. 
Under this framework, we obtain the first preliminary result towards the construction of proper test procedures.

\begin{theorem}\label{Theorem1}
Suppose Assumption (A1) holds. Then, as $N\to \infty$, we have convergence in distribution
\[\sqrt{N} \vC(\widehat \vv -\vv)\stackrel{\mathcal {D}}{\longrightarrow} {\mathcal{N}_{m}\left(\vnull_{m},\vC\vSigma \vC^\top\right)},\]
where  $\vSigma=\bigoplus_{i=1}^a {\kappa_i}^{-1}\cdot  \vSigma_i$ and $\vSigma_i=\Cov(\vech(\vepsilon_{i1}\vepsilon_{i1}^\top)) $ for $i=1,\dots, a$.
\end{theorem}
{Together with a consistent estimator for (all or certain parts of) $\vSigma$, this result will allow us to develop asymptotic tests for the null hypothesis \eqref{eq:hypo cov}}. { 
To this end, we define the empirical estimator $\widehat \vSigma: =\bigoplus_{i=1}^a {N}/{n_i}  \cdot \widehat \vSigma_i$ for $\vSigma$, where
$$\widehat \vSigma_i=\frac{1}{n_i-1} \sum\limits_{k=1}^{n_i}\left[\vech\left(\widetilde \vX_{ik}\widetilde \vX_{ik}^\top- \sum\limits_{\ell=1}^{n_i}\frac{\widetilde \vX_{i\ell}\widetilde \vX_{i\ell}^\top}{n_i}\right)\right]\left[\vech\left(\widetilde \vX_{ik}\widetilde \vX_{ik}^\top- \sum\limits_{\ell=1}^{n_i}\frac{\widetilde \vX_{i\ell}\widetilde \vX_{i\ell}^\top}{n_i}\right)\right]^\top.
$$
Here, $\widetilde \vX_{ik}:=\vX_{ik}-\overline \vX_{i\cdot}$ denotes the centered version of observation $k$ in group $i$. The consistency of the matrices $\widehat \vSigma_i$ for $\vSigma_i$ and thus of $\widehat \vSigma$ is established in the supplementary material.}

{Now potential test statistics may lean on well-known quadratic forms used for mean-based MANOVA-analyses  in heteroscedastic designs (\cite{kon:2015, bathke2018}). To unify several approaches we consider 
\begin{equation}\label{eq:general_QF}
\widehat Q_{\vv}=N\left[ \vC\widehat \vv - \vzeta\right]^\top \vE(\vC,\widehat{\vSigma})\left[  \vC\widehat \vv - 
\vzeta\right],
\end{equation}
where, $\vE(\vC,\widehat{\vSigma}) \in \R^{m\times m}$ is some symmetric matrix that can be written as a function of the hypothesis matrix $\vC\in \R^{m \times ap}$ and the covariance matrix estimator $\widehat{\vSigma}\in\R^{ad\times ad}$. 
{In order to analyze the limit behaviour of $\widehat Q_{\vv}$ we assume throughout that $\vE(\vC,\widehat \vSigma)\stackrel{\mathcal{P}}{\to}\vE(\vC,\vSigma)$ holds which is, e.g., fulfilled if $\vE$ is continuous in its second argument. }
{Choices covered by this general formulation include the following:}
\begin{itemize}
\item[1.] An {\bf ANOVA-type-statistic (ATS):}~
$ATS_{\vv}(\widehat{\vSigma})= N\left[  \vC\widehat \vv - \vzeta\right]^\top\left[  \vC\widehat \vv - 
\vzeta\right]/\tr\left(\vC\widehat{\vSigma}\vC^\top\right)$ corresponding to 
$\vE(\vC,\widehat \vSigma) = \vI_m/\tr(\vC \widehat \vSigma \vC^\top)$.
\item[2.] A {\bf Wald-type-statistic (WTS):}~ 
$WTS_{\vv}(\widehat{\vSigma}) = N\left[  \vC\widehat \vv - \vzeta\right]^\top\left(\vC \widehat{\vSigma} \vC^\top\right)^+\left[  \vC\widehat \vv - \vzeta\right]$. Here, $\vE(\vC,\widehat{\vSigma}) = \left(\vC  \widehat\vSigma \vC^\top\right)^+$ is the Moore-Penrose-inverse of $\vC  \widehat\vSigma \vC^\top$. As we will see later, the usual $\chi_{f}^2$-limit distribution with $f=\rank(\vC)$ will appear under the additional assumption $\vSigma>0$. This or comparable conditions are required to garantee $\vE(\vC,\widehat \vSigma)\stackrel{\mathcal{P}}{\to}\vE(\vC,\vSigma)$.
\item[3.] Substituting $\widehat\vSigma$ in the WTS with $\widehat\vSigma_0$, the 
diagonal matrix only containing the diagonal elements of $\widehat\vSigma$, leads to the so-called {\bf modified ANOVA-type statistic (MATS)} given by $MATS_{\vv}(\widehat\vSigma) = N\left[  \vC\widehat \vv - \vzeta\right]^\top\left(\vC \widehat \vSigma_0 \vC^\top\right)^+\left[  \vC\widehat \vv - \vzeta\right]$. To study its asymptotics we need to assume $\vSigma_0>0$.
\end{itemize}
In 2. and 3. the additional assumptions are needed to ensure that the inner Moore Penrose inverse is consistent. The following result establishes the asymptotic distribution of all quadratic forms of type \eqref{eq:general_QF} and covers all the cases 1.-3..}
\begin{theorem}\label{Verteilung}  
{Under Assumption (A1) and the null hypothesis $\mathcal{H}_0^v:\vC\vv=\vzeta$, the quadratic form $\widehat Q_{\vv}$ defined by \eqref{eq:general_QF}
 has, asymptotically, a ``weighted $\chi^2$-distribution''. That is, 
$$\widehat Q_{\vv}\stackrel{\mathcal {D}}{\To} \sum_{\ell=1}^{ap} \lambda_\ell B_\ell,
$$
where $B_\ell \stackrel{i.i.d.} {\sim} \chi_1^2$ and $\lambda_\ell, \ell =1,\dots, ap,$ are the eigenvalues of $(\vSigma^{1/2}\vC^\top\vE(\vC,\vSigma)\vC \vSigma^{1/2})$.}
\end{theorem}

This result allows the definition of a natural test procedure in the WTS given by {$\varphi_{WTS}=\ind\{ {WTS_{\vv}(\widehat \vSigma)} \notin (-\infty,\chi^2_{f;1-\alpha}]\}$}. However, the additional condition \textcolor{black}{$\vSigma>0$}, ensuring asymptotic correctness of $\varphi_{WTS}$, may not always be satisfied in practice.

Since this condition is not needed for the ANOVA-type statistic $A_N = ATS_{\vv}(\widehat \vSigma)$, we focus on the ATS in what follows; noting that the MATS did also show good finite sample properties in simulations, see the supplement for details. 
As the limit distribution of the ATS depends on unknown quantities, we cannot calculate critical values from \Cref{Verteilung} directly. 
To this end, we employ resampling techniques for calculating proper critical values. We thereby focus on two resampling procedures: a parametric and a wild bootstrap as both methods have shown favorable finite sample properties in multivariate mean-based MANOVA %
(\cite{kon:2015}, \citet*{friedrich2016}, \cite{friedrich2017mats}, and \citet*{zimmermann2019}).
That these procedures also lead to valid testing procedures in the current setting is proven in the subsequent section.

\section{\textsc{Resampling Procedures}}\label{Resampling Procedures}

To derive critical values for the non-pivotal test statistics like $ATS_{\vv}$, we consider two common kinds of bootstrap techniques: a parametric and a wild bootstrap as applied for heteroscedastic MANOVA. Since we deal with covariances instead of expectations, some adjustments have to be made in order to prove their asymptotic correctness.

\subsection{Parametric Bootstrap}\label{Parametric Bootstrap}
To motivate our first resampling strategy, note that  
\bqan\label{eq:parametric bootstrap}
{\sqrt{N}}(\widehat \vv_i  -\vv_i)= {\sqrt{N}}\vech\left(\frac 1 {n_i-1} \sum\limits_{k=1}^{n_i}\left[\vep_{i k}\vep_{i k}^\top-{\vV}_i\right]\right)+\lan_P(1)\stackrel{\mathcal{ D}}{\to} \mathcal{N}_p\left(\vnull_{p},\frac 1 {\kappa_i}\vSigma_i\right)\eqan
 follows from the proof of \Cref{Theorem1}. 

Thus, to mimick its limit distribution and afterwards the structure of the test statistic, we generate bootstrap vectors $\vY_{i 1}^*,...,\vY_{i n_i}^*\stackrel{i.i.d.}{\sim} \mathcal N_{p}\left(\vnull_p,\widehat{\vSigma}_{i}\right),$ for given realisations $\vX_{i 1},..., \vX_{i n_i}$ with estimators   $\widehat \vSigma_{i}$. We then calculate 
{$\widehat \vSigma_i^*$, the empirical covariance matrix of the bootstrap sample $\vY_{i1}^*,...,\vY_{in_i}^*$   and set $\widehat \vSigma^*:=\bigoplus_{i=1}^a{N}/{n_i}\cdot  \widehat \vSigma_i^*$.

} The next theorem ensures the asymptotic correctness of this approach.
\begin{theorem}\label{PBTheorem1}Under Assumption (A1), the following results hold:\\
(a) For $i=1,...,a$, the conditional distribution of $\sqrt{N}\ \overline \vY_i^*$, given the data, converges weakly to $ \mathcal{N}_{p}\left(\vnull_p,{\kappa_i}^{-1}\cdot \vSigma_i\right)$ in probability. {Since $\widehat \vSigma_i^*\to \vSigma_i$ in probability, the unknown covariance matrix $\vSigma_i$ can be estimated through $\widehat \vSigma_i^*$.}\\
(b) The conditional distribution of $\sqrt{N}\ \overline \vY^*$, given the data, converges weakly to \\$ \mathcal{N}_{a\cdot p}\left(\vnull_{a\cdot p},\bigoplus_{i=1}^a {\kappa_i}^{-1}\cdot \vSigma_i\right)$ in probability. { Since $\widehat \vSigma^*\to \vSigma$ in probability, the unknown covariance matrix $\vSigma$ can be estimated through $\widehat \vSigma^*$.}

\end{theorem}

{As a consequence, it is reasonable to calculate the bootstrap version of the general quadratic form \eqref{eq:general_QF}} as
$Q_{\vv}^*=N[ \vC\overline \vY^*  ]^\top \vE(\vC,\widehat \vSigma^*)[ \vC\ \overline \vY^*  ].$
For the ATS, e.g., this leads to $ATS_{\vv}^*=N[ \vC\overline \vY^*  ]^\top[ \vC\ \overline \vY^*  ]\big/\tr(\vC \widehat \vSigma^* \vC^\top).$
{The bootstrap versions approximate the null distribution of $\widehat Q_{\vv}$, as established below.}
\begin{corollary}\label{KorParametric}
For each parameter $\vv\in \R^{a\cdot p}$ and $\vv_0$ with $\vC\vv_0=\vzeta$, we have under Assumption (A1) that
\[\begin{array}{l}
\sup\limits_{x\in\R}\big\lvert P_{\vv}(Q_{\vv}^*\leq x\lvert \vX)-P_{\vv_0}(\widehat Q_{\vv}\leq x)\big\lvert \stackrel{\mathcal P}{\to}0,\\[1.5ex]
\end{array}\] where $P_{\vv}$ denotes the (un)conditional distribution of the test statistic when $\vv$ is the true underlying  vector.
 \end{corollary}

Denoting with $c_{ATS^*, 1-\alpha}$  the $(1-\alpha)$-quantile of the conditional distribution of $ATS_{\vv}^*$ given the data, we obtain 
$\varphi_{ATS}^*=\ind\{ ATS_{\vv}(\widehat \vSigma) \notin (-\infty,c_{ATS^*,1-\alpha}]\}$ as asymptotic level $\alpha$ test.
 
Beyond being helpful to carry out an asymptotic level $\alpha$ test in the $ATS_{\vv}$, resampling can also be used to enhance the finite sample properties of the $WTS_{\vv}$. In fact, utilizing \Cref{PBTheorem1}  shows that a parametric bootstrap version of the $WTS_{\vv}$, say $WTS_{\vv}^*$, is also asymptotically $\chi_{\rank(\vC)}^2$-distributed, under the assumption given in \Cref{Theorem1}.
Thus, it leads to a valid parametric bootstrap $WTS_{\vv}$-test as long as $\vSigma_i>0$ for all $i=1,\dots,a$.

\subsection{Wild Bootstrap}
As a second resampling approach, we consider the wild bootstrap. Hereby the structure of the data is kept more than for the parametric bootstrap since no fixed distribution is used.%
 In the mean-based analysis, convenient wild bootstrap multipliers are multiplied with the realizations to get the bootstrap sample. In contrast, we have to multiply them with $p$-dimensional random vectors of the kind $\vech(\vX_{ik}\vX_{ik}^\top)$, to ensure asymptotic correctness due to \eqref{eq:parametric bootstrap}.

Specifically, generate i.i.d. random weights  $W_{i1},...,W_{i n_i}$, $i=1,...,a,$ independent of the data, with $\E(W_{i1})=0$ and $\Var(W_{i1})=1$. Common choices are for example standard  distributed random variables or random signs. 
Afterwards the wild bootstrap sample is defined as $\vY_{i k}^\star=W_{ik}\cdot  \left[\vech(\widetilde \vX_{ik} \widetilde \vX_{ik}^\top)-n_i^{-1} \sum_{\ell=1}^{n_i}\vech(\widetilde \vX_{i\ell}\widetilde {\vX}_{i\ell}^\top) \right]$, where again centering is needed to capture the correct limit structure.
Defining $\widehat{\vSigma}_i^\star$ as the empirical covariance matrix of $\vY_{i1}^\star,...,\vY_{i n_i}^\star$ and setting $\widehat{\vSigma}^\star=\bigoplus_{i=1}^a N/n_i\cdot \widehat{\vSigma}_i^\star$, we obtain the following theorem.

\begin{theorem}\label{WBTheorem}
Under Assumption (A1), the following results hold:\\
(a) For $i=1,...,a$, the conditional distribution of $\sqrt{N} \ \overline \vY_i^\star$, given the data converges weakly to $ \mathcal{N}_{p}\left(\vnull_p,{\kappa_i}^{-1}\cdot  \vSigma_i\right)$ in probability. 
{Since $\widehat \vSigma_i^\star \to \vSigma_i$ in probability, the unknown covariance matrix $\vSigma_i$ can be estimated through $\widehat \vSigma_i^\star$.}
\\
(b) The conditional distribution of $\sqrt{N}\ \overline \vY^\star$, given the data converges weakly to \\$ \mathcal{N}_{a\cdot p}\left(\vnull_{a\cdot p},\bigoplus_{i=1}^a  {\kappa_i}^{-1} \cdot  \vSigma_i\right)$ in probability. { Since $\widehat \vSigma^\star \to \vSigma$ in probability, the unknown covariance matrix  $\vSigma$ can be estimated through $\widehat \vSigma^\star$. }
\end{theorem}

The result again gives rise to define a wild bootstrap quadratic form \\
$Q_{\vv}^\star=N[ \vC\overline \vY^\star  ]^\top \vE(\vC,\widehat \vSigma^\star)[ \vC\ \overline \vY^\star  ], $ 
where, e.g., an $ ATS_{\vv}(\widehat \vSigma)$ wild bootstrap counterpart is given by
$ATS_{\vv}^\star=N[ \vC\overline \vY^\star ]^\top[ \vC \overline \vY^\star  ]\big/\tr(\vC \widehat \vSigma^\star \vC^\top).
$ Similar to the parametric bootstrap, the next theorem guarantees the approximation of the original test statistic by its bootstrap version.

\begin{corollary}\label{KorWild}
Under the assumptions of \Cref{KorParametric}, we have convergence
\[\begin{array}{l}
\sup\limits_{x\in\R}\textcolor{black}{\big\lvert P_{\vv}(Q_{\vv}^\star\leq x\lvert \vX)-P_{\vv_0}(\widehat Q_{\vv}\leq x)\big\lvert \stackrel{\mathcal P}{\to}0.}
\end{array}\] 
 \end{corollary}
Therefore, analogous to $\varphi_{ATS}^*$, we define $\varphi_{ATS}^\star:=\ind\{ ATS_{\vv}(\widehat \vSigma) \notin (-\infty,c_{ATS^\star,1-\alpha}]\}$  as asymptotic level $\alpha$ test, with $c_{ATS^\star, 1-\alpha}$ denoting the $(1-\alpha)$ quantile of the conditional distribution of $ATS_{\vv}^\star$ given the data. 

Similar wild bootstrap versions of the $WTS_{\vv}$ or comparable statistics can again be defined and used to calculate critical values if $\vSigma>0$ is fulfilled, see Section~\ref{Simulations} below for the WTS and the supplement for another, less known, possibility.

\section{\textsc{Simulations}}\label{Simulations}
The above results are valid for large sample sizes. For an evaluation of the finite sample behavior of all methods introduced above, we have conducted extensive simulations regarding
\begin{itemize}
\item[(i)] their ability to keep the nominal significance level and
\item[(ii)]their power to detect certain alternatives in various scenarios.
\end{itemize}

In particular, we studied three different kinds of hypotheses:
\begin{itemize}
\item[$A$)] {Equal Covariance Matrices:\ } $\mathcal{H}_0^{\vv}:\vV_{1}=\vV_{2}$ with $a=2$ groups. 
\item[$B$)] Equal Diagonal Elements: \ $\mathcal{H}_0^{\vv}:\vV_{1 11}=...=\vV_{1 dd}$ in the one sample case.

\item[$C$)] Trace Test: \ $\mathcal{H}_0^{\vv}:\tr(\vV_{1})=\tr(\vV_{2})$ with $a=2$ groups.

\end{itemize}

Each of these hypotheses can be formulated with a proper projection matrix $\vC$. While $\vC(A)=\vP_2\otimes \vI_d$ and $\vC(C)= \vP_2\otimes [\vh_d\cdot \vh_d^\top]/d$ follows directly from Section 2, $\vC(B)=\diag(\vh_d)-\vh_d\cdot \vh_d^\top/d$ is an adaptation of  $\vP_d$.\\
For each hypothesis, we have simulated the two bootstrap methods based on the ANOVA-type statistic $\varphi_{ATS}^*$ and $\varphi_{ATS}^\star$, as well as the Wald-type-statistic $\varphi_{WTS}^*$ {and} $ \varphi_{WTS}^\star$. The latter ones are based on the parametric bootstrap version of the WTS, given by
\begin{equation}\label{WTSP} WTS^*(\widehat \vSigma^*):=N \left[\vC\overline \vY^*\right]^\top \left(\vC \widehat \vSigma^*\vC^\top\right)^+ \left[\vC\overline \vY^*\right]\end{equation}
and the  wild bootstrap version given by
\begin{equation}\label{WTSW} WTS^\star(\widehat \vSigma^\star):=N \left[\vC\overline \vY^\star\right]^\top \left(\vC \widehat \vSigma^\star\vC^\top\right)^+ \left[\vC\overline \vY^\star\right].\end{equation}

Moreover, the asymptotic version $\varphi_{WTS}$ based upon the $\chi_{\rank(\vC)}^2$-approxi\-mation serves as another competitor. 
{As additional competitor, we consider a Monte-Carlo test in the ATS. Recall that its limiting null distribution is given by $A_0:=\sum_{k=1}^{m} \lambda_k B_k\big/\tr\left(\vC \vSigma \vC^\top\right)$ for $B_k\stackrel{i.i.d.}{\sim}\chi_1^2$ and $\lambda_k\in eigen\left(\vC\vSigma \vC^\top\right)$. Plugging in $\widehat{\vSigma}$ for $\vSigma$ 
and repeatedly generating $C_k$'s within 10.000 Monte-Carlo of $A_0$, we obtain an estimated 
$(1-\alpha)$-quantile $q_{1-\alpha}^{MC}$ of the distribution of $A_0$. This finally defines the Monte-Carlo ATS test $\varphi_{ATS}:=\ind\{ ATS_{\vv}(\widehat \vSigma) \notin (-\infty,q_{1-\alpha}^{MC}]\}$.

In the special case of scenario $A)$, we have also considered the tests from \\\cite{zhang1992,zhang1993} based on Bartlett's test statistic, 
along with a so-called separate bootstrap as well as a pooled bootstrap to calculate critical values. 
We denote these tests by $\varphi_{B-S}$ and $\varphi_{B-P}$.
While the first is asymptotically valid under the same conditions as our tests, the pooled bootstrap procedure additionally requires 
$\E\left(\left[vech(\vepsilon_{1}\vepsilon_{1}^\top)\right]\left[vech(\vepsilon_{1}\vepsilon_{1})^\top\right]^\top  \right)=\E\left(\left[vech(\vepsilon_{2}\vepsilon_{2}^\top)\right]\left[vech(\vepsilon_{2}\vepsilon_{2}^\top)\right]^\top  \right)$. \\

Additionally, we simulated Box's M-test as it is the most popular test for scenario $A)$, although it requires normally distributed data. 
There are two common ways to determine critical values for this test (\cite{box1949}): Utilizing a $\chi_f^2$-approximation with $f=\rank(\vC)$ degrees of freedom or an $F$-approximation with estimated degrees of freedom. For ease of completeness, we decided to simulate both.
\\
 On an abstract level, the hypotheses considered thus far also fall into the framework presented by  \cite{zhang1993}. However, they do not provide concrete test statistics that we could use for comparison purposes. Other existing tests, such as the one by \cite{gupta2006} rely on rather different model assumptions, which also makes a comparative evaluation difficult.
All simulations were conducted by means of the \textsc{R}-computing environment version 3.6.1 \cite{R} 
with $N_{sim}=2\cdot 10^4$ runs, 1000 bootstrap runs and $\alpha=5\%$.
\subsection*{Data generation}
We considered $5$-dimensional observations generated independently according to the model $\vX_{ik} = \mu_i + \vV^{1/2} \vZ_{ik}, i=1,\dots, a, k=1,\dots,n_i$ with $\vmu_1=(1^2,2^2,...,5^2)/4$ and $\vmu_2=\vnull_5$. Here, the marginals of $\vZ_{ik} = (Z_{ikj})_{j=1}^5$ were either simulated independently from 
\begin{itemize}

\item a standard normal distribution, i.e. $Z_{ikj} \sim  \mathcal {N}(0,1)$
\item a standardized centered gamma distribution i.e. $(\sqrt{2} Z_{ikj}+2)\sim \mathcal G(2,1)$
\item a standardized centered skew normal distribution with location parameter $\xi=0$, scale parameter $\omega=1$ and $\alpha=4$. The density of a skew normal distribution is given through $\frac{2}{\omega} \phi\left(\frac{x-\xi}{\omega}\right)\Phi\left(\alpha\left(\frac{x-\xi}{\omega}\right)\right)$, where $\phi$ denotes the densitiy of a standard normal distribution and $\Phi$ the according distribution function.
\end{itemize}
For the covariance matrix, an autoregressive structure with parameter $0.6$  was chosen, i.e., $(\vV)_{ij}=0.6^{|i-j|}$. More simulation results with different covariance matrices {and more distributions} can be found in the supplement. 
This includes hypotheses with more groups or settings with a higher dimension of the observations.

Note that the chosen dimension of $d=5$ leads to an effective dimension of $p=15$ of the unknown parameter (i.e., covariance matrix) in each group. Hence in scenario $A)$, the vector $\vnu$ defining the null hypothesis \eqref{eq:hypo cov} actually consists of $30$ unknown parameters. 
To address this quite large dimension, we considered four different small to large total sample sizes of {$N\in\{50, 100, 250, 500\}$}. Moreover, in scenario $A)$ and $C)$ these were divided into two groups by setting 
$n_1=0.6\cdot N$ and $n_2=0.4\cdot N$. In scenario $B)$ the sample size is {$n\in \{25,50,125,250\}$}
Thus, we had between 20 and 300 independent observations to estimate the unknown covariance matrix in each group.

\subsection{Type-I-error }
The following tables display the simulated type-I-error rates for all these settings. Values inside the $95\%$ binomial interval $[0.047; 0.053]$ are printed bold.

\begin{table}[h!]
      \centering
      \begin{scriptsize}

    \begin{tabular}{|l|| c|c|c|c||c|c|c|c||c|c|c|c|}\hline
   
      &\multicolumn{4}{|c||}{Normal}&\multicolumn{4}{|c||}{Skewed Normal}&\multicolumn{4}{|c|}{Gamma}\\\cline{1-13}   
  \hspace*{.1cm}N  &   50&100&250&500&   50&100&250&500&  50&100&250&500
       \\\hline\hline
\hspace{-0.1cm}ATS-Para &.0579 & { .0540} & {\bf .0518} & {\bf .0515} & .0589 & { .0538} &{\bf .0528} &{\bf .0488} &{\bf .0485} & .0439 & .0439 &{ .0464} \\    \hline
\hspace{-0.1cm}ATS-Wild &.0797 & .0672 &{ .0558} &{ .0533} & .0915 & .0708 & .0619 & .{\bf 0.522} & .0995 & .0784 & .0611 &{ .0552} \\    \hline
\hspace{-0.1cm}ATS & .0634 & .0562 &{\bf .0520} &{\bf .0510} & .0640 &{ .0543} &{\bf .0530} &{\bf .0484} &{ .0538} &{ .0462} &{ .0447} & { .0458} \\    \hline
\hspace{-0.1cm}WTS-Para & .0659 & .0661 & .0623 & .0566 & .0798 & .0727 & .0648 & .0604 & .0800 & .0690 & .0638 & .0582 \\    \hline
\hspace{-0.1cm}WTS-Wild  & .0961 & .0852 & .0706 & .0612 & .1167 & .0975 & .0786 & .0689 & .1300 & .1083 & .0870 & .0707 \\    \hline
\hspace{-0.1cm}WTS-$\chi_{15}^2$  & .5000 & .2161 & .1054 & .0757 & .5231 & .2387 & .1100 & .0812 & .5448 & .2389 & .1085 & .0764 \\    \hline
\hspace{-0.1cm}Bartlett-S  &.0111 & .0371 &{\bf .0478} &{\bf .0485} & .0166 & .0400 &{\bf .0528} &{\bf .0515} & .0264 & .0594 & .0655 & .0613 \\    \hline
\hspace{-0.1cm}Bartlett-P & .0199 & .0360 &{ .0452} &{ .0467} & .0254 & .0361 &{ .0452} &{\bf .0480} & .0299 & .0405 &{ .0451} &{\bf .0485} \\    \hline
\hspace{-0.1cm}Box's M-$\chi_{15}^2$\hspace{-0.15cm} & .0638 & .0575 &{\bf .0521} &{\bf .0496} & .1075 & .0976 & .0956 & .0938 & .2707 & .2896 & .3156 & .3250 \\    \hline
\hspace{-0.1cm}Box's M-F& .0609 & .0567 & {\bf .0520} &{\bf .0496} & .1012 & .0961 & .0952 & .0938 & .2612 & .2881 & .3153 & .3249 \\    \hline

\end{tabular}
  \caption{Simulated type-I-error rates ($\alpha=5\%$) in scenario $A)$ ($\mathcal{H}_0^{\vv}:\vV_{1}=\vV_{2}$)  for ATS, WTS, MATS, Bartlett's test and Box's M-test, always with the same relation between group sample sizes by $n_1:=0.6\cdot N$ resp. $n_2:=0.4\cdot N$. The 5-dimensional observation vectors have the covariance matrix $(\vV)_{ij}=0.6^{|i-j|}$.}
   \label{tab:SimAH}

\end{scriptsize}
  
    \end{table}
 
In almost all simulation settings, the wild bootstrap led to more liberal results, whereas the parametric bootstrap was also liberal for the WTS but had no clear tendency for the ATS. 
For larger sample sizes, the ATS with critical values based on the weighted sum of $\chi^2$ random variables behaved similarly to the ATS with parametric bootstrap, 
while for smaller sample sizes, the simulated type-I-error rates differed more from the nominal $\alpha$-level.  \\

Overall, the results of the ATS were preferable compared to the WTS. 
This matches the conventional wisdom that the WTS generally exhibits a liberal behavior and requires large sample sizes to perform well. Moreover, the WTS requires the condition on the rank of $\vSigma$, which is difficult to check in practice because of the special structure of $\vSigma$. 
In contrast, the ATS is capable to handle all these scenarios.

Therefore, it remains to compare these tests with those based on Bartlett's statistic.\\
The additional condition required for the pooled bootstrap is fulfilled. Therefore, \Cref{tab:SimAH} contains also the results of $\varphi_{B-P}$.

\begin{table}[tbp]      \begin{scriptsize}
    \centering
\begin{tabular}{|l||c|c|c|c||c|c|c|c||c|c|c|c|}\hline
    &\multicolumn{4}{|c||}{Normal}&\multicolumn{4}{|c||}{Skewed Normal}&\multicolumn{4}{|c|}{Gamma}\\\cline{1-13} 
  \hspace*{.1cm}N  &    25&50&125&250&  25&50&125&250&  25&50&125&250 
       \\\hline\hline

\hspace{-0.1cm}ATS-Para  & {.0465} & {\bf.0473} & {\bf.0495} &{\bf .0505} &{\bf .0481} & .0419 &{ .0454} & {\bf.0483} & .0388 & .0363 & .0371 & .0407 \\    \hline
\hspace{-0.1cm}ATS-Wild& .0682 & .0573 & {.0542} & {\bf.0527} & .0787 & .0618 & {.0550} &{ .0547} & .0805 & .0645 & {.0535} &{\bf .0524} \\    \hline
  \hspace{-0.1cm}ATS & {.0547} & {\bf.0501} & {\bf.0492} & {\bf.0501} & .0566 & {.0451} & {.0458} & {\bf.0487} & {.0455} & .0383 & .0373 & .0397 \\    \hline
\hspace{-0.1cm}WTS-Para & .0855 & .0702 & .0622 & {.0545} & .1099 & .0886 & .0726 & .0618 & .1441 & .1136 & .0839 & .0711 \\    \hline
\hspace{-0.1cm}WTS-Wild & .1052 & .0795 & .0660 & {.0557} & .1458 & .1112 & .0826 & .0675 & .2099 & .1590 & .1076 & .0847 \\    \hline
  \hspace{-0.1cm}WTS-$\chi_{4}^2$ &  .1826 & .1109 & .0682 & .0594 & .2207 & .1277 & .0797 & .0708 & .2609 & .1628 & .0939 & .0761 \\     \hline

\end{tabular}
  \caption{Simulated type-I-error rates ($\alpha=5\%$) in scenario $B)$ ($\mathcal{H}_0^{\vv}:\vV_{1 11}=...=\vV_{1 55}$) for ATS and WTS. The 5-dimensional observation vectors have the covariance matrix $(\vV)_{ij}=0.6^{|i-j|}$.  }
   \label{tab:SimBH}
     \end{scriptsize}
    \end{table}

\begin{table}[tbp]
   \begin{scriptsize}
    \centering
  \begin{tabular}{|l||c|c|c|c||c|c|c|c||c|c|c|c|}\hline
      &\multicolumn{4}{|c||}{Normal}&\multicolumn{4}{|c||}{Skewed Normal}&\multicolumn{4}{|c|}{Gamma}\\\cline{1-13}   
     \hspace*{.1cm}N    & 50&100&250&500&   50&100&250&500&   50&100&250&500
       \\\hline\hline

\hspace{-0.1cm}ATS-Para & .0651 & .0581 &{ .0537} &{ .0539} & .0690 & .0589 & {\bf.0530} & {\bf.0514} & .0715 & .0628 & {.0540} & {.0552} \\ \hline
  \hspace{-0.1cm}ATS-Wild& .0686 & .0598 & {.0542} & {.0545} & .0738 & .0621 & {.0540} & {\bf.0521} & .0848 & .0688 & {.0550} & {.0556} \\ \hline
  \hspace{-0.1cm}ATS & .0739 & .0609 & {.0544} & {.0541} & .0779 & .0623 & {.0540} & {\bf.0517} & .0814 & .0655 &{ .0540} & {.0538} \\ \hline
  \hspace{-0.1cm}WTS-Para  & .0651 & .0581 & {.0537} & {.0539} & .0690 & .0589 & {\bf.0530} &{\bf .0514} & .0715 & .0628 & {.0540} & {.0552} \\ \hline
  \hspace{-0.1cm}WTS-Wild  & .0686 & .0598 & {.0542} & {.0545} & .0738 & .0621 & {.0540} & {\bf.0521} & .0848 & .0688 & {.0550} & {.0556} \\ \hline
 \hspace{-0.1cm}WTS-$\chi_{1}^2$& .0736 & .0605 & {.0535} & {.0538} & .0775 & .0619 & {.0538} & {\bf.0518} & .0811 & .0651 & {.0540} & {.0540} \\ \hline

\end{tabular}
  \caption{Simulated type-I-error rates ($\alpha=5\%$) in scenario $C)$ ($\mathcal{H}_0^{\vv}:\tr(\vV_{1})=\tr(\vV_{2})$)  for ATS, WTS, MATS, Bartlett's test and Box's M-test, always with the same relation between group sample sizes by $n_1:=0.6\cdot N$ resp. $n_2:=0.4\cdot N$. The 5-dimensional observation vectors have the covariance matrix $(\vV)_{ij}=0.6^{|i-j|}$.}
  \label{tab:SimCH}
     \end{scriptsize}
\end{table}
For all distributions, $\varphi_{ATS}^*$ showed good results especially for normal distribution and skewed normal distribution where the type-I-error rate was always better than those of 
$\varphi_{B-S}$ and $\varphi_{B-S}$. Also for the gamma distribution $\varphi_{B-S}$ performed worse while for  bigger N  the simulated error-rates  of $\varphi_{ATS}^*$  and $\varphi_{B-S}$ were comparable. Most of all  $\varphi_{ATS}^*$ provided good values for small samples, while both tests based on a Bartlett statistic needed large sample sizes.
At last, the popular Box's M-test worked quite well under normality but showed poor results (type-I-error rates of more than 20\%) when this condition was violated. This sensitivity to the violation of normal distribution may have the consequence in practice that small p-values could be untrustworthy, independent of whether $\chi^2$ or F distribution was used. 
But also for normality, the performance was not essentially better than $\varphi_{ATS}$ and (with small exceptions) clearly worse than $\varphi_{ATS}^*$.
This also underlines the benefit of the newly proposed test for this popular null hypothesis.

Moreover, the resampling procedure used in \cite{zhang1993}  occasionally encountered covariance matrices without full rank, especially for smaller sample sizes. This creates issues in the algorithm because the determinant of these matrices is zero, and the logarithm at this point is not defined. Regretfully this situation wasn't discussed in the original paper, so we just excluded these values. Certainly, this would constitute a drastic user intervention in applying the bootstrap and also influencing the conditional distribution. Nevertheless, it was necessary to use this adaptation in all our simulations containing these tests.
This effect can also occur in Box's M-test, but comparatively rarely because there is no bootstrap involved.
\\\\
All in all, in scenario $A)$ the $ATS^*$ and the Monte-Carlo ATS test exhibited the best performance over all distributions and, in particular small sample sizes.

For scenario $B)$ the results in \Cref{tab:SimBH} again show the good performance of $\varphi_{ATS}^*$ for small sample sizes. With the exception of the gamma distribution, where for large sample sizes  $\varphi_{ATS}^\star$ had an error rate closer to our $\alpha$ level, the ATS using the parametric bootstrap approach had by far the best results.

At last \Cref{tab:SimCH} shows the results from scenario $C)$. Due to the fact that the rank of the hypothesis matrix is 1, there is no difference between the WTS and the ATS. All our tests $\varphi_{ATS}^*,\varphi_{ATS}^\star$ and $\varphi_{ATS}$ showed comparable results while again $\varphi_{ATS}^*$ had the best small sample performance. In comparison to the other scenarios, the error rates were a bit worse than before. But we have to take into account that this is the most challenging hypothesis, which only considers the diagonal elements of the covariance matrix. Nevertheless, the results for sample sizes 250 and 500 were convincing.

The effect of using other types of covariance matrices, which is considered in the supplement, was not significant and not systematic. Therein, we also investigated testing for a given covariance matrix.  
 Here, only the type-I-error rate of the ANOVA-type statistic with critical values obtained from the parametric bootstrap and the Monte-Carlo ATS test showed sufficiently good results.\\
To sum up, we only recommend the use of any of the three tests based on the ATS. 
All three exhibited good simulation results for comparably small sample sizes and are (asymptotically) valid without additional requirements on $\vSigma$.  
Additional simulations, given in the supplementary material, also confirm this conclusion, especially for a higher dimension or more groups.

\subsection{Power}
For a power simulation, it is unfortunately not possible to merely shift the observations by a proper vector to control the distance from the null hypothesis.  
Thereto we have multiplied the observation vectors $\vX$ with a proper diagonal matrix, given by  $ \vDelta=\vI_d+\diag(1,0,...,0)\cdot \delta$ for $\delta\in[0,3]$.
This was associated with a one-point-alternative that is known from testing expectation vectors to be challenging, namely a deviation in just one component, which is usually difficult to detect.
 
 \begin{figure}[H]
\includegraphics[trim= 7mm 13mm 5mm 20mm,clip,scale=0.735]{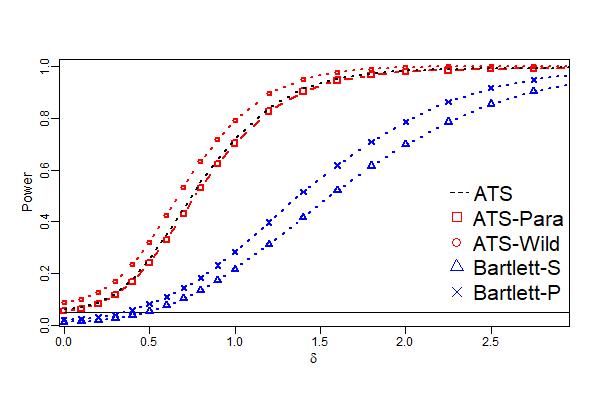}
\caption{Simulated power in scenario $A)$ ($\mathcal{H}_0^{\vv}:\vV_{1}=\vV_{2}$) for {ATS with wild bootstrap, parametric bootstrap, and Monte-Carlo critical values,} as well as the test based on Bartlett's statistic with separate and pooled bootstrap.
The  5-dimensional  vectors were based  on the skewed normal distribution, with covariance matrix  $(\vV)_{ij}=0.6^{|i-j|}$ and $n_1=30, n_2=20$. The considered alternative is a one-point-alternative.}
\end{figure}
In this way  $\vC\vech\left(\vDelta \vV \vDelta^\top\right)-\vzeta\neq \vnull$, were $n_1+n_2=50$ was used  to investigate small size behavior, while the dimension was again $d=5$, leading to $p=15$. Moreover for a second alternative the observation vectors $\vX$  were multiplied by $ \vDelta=\vI_d+\diag(1,2,...,d)/d\cdot \delta$ for $\delta\in[0,3]$, which corresponds to a so-called trend-alternative. Due to computational reasons and because of the performance under the null hypothesis described in the last section, we have only investigated the power of $\varphi_{ATS}^*$, $\varphi_{ATS}^\star$ and $\varphi_{ATS}$ as well as $\varphi_{B-P}$ and $\varphi_{B-S}$  from \cite{zhang1993} for skewed normal distributed random variables.

\begin{figure}[H]
\includegraphics[trim= 7mm 13mm 5mm 20.8mm,clip,scale=0.735]{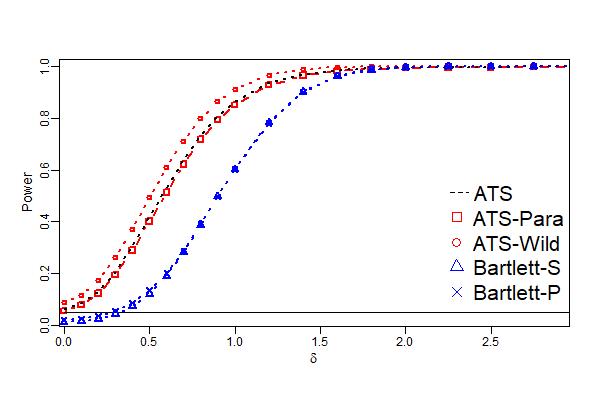}
\caption{Simulated power in scenario $A)$ ($\mathcal{H}_0^{\vv}:\vV_{1}=\vV_{2}$) for {ATS with wild bootstrap, parametric bootstrap, and Monte-Carlo critical values,} as well as the tests based on Bartlett's statistic with separate and pooled bootstrap.
The  5-dimensional  vectors were based  on the skewed normal distribution, with covariance matrix  $(\vV)_{ij}=0.6^{|i-j|}$ and $n_1=30, n_2=20$. The considered alternative is a trend-alternative.}
\end{figure}
Overall, the ATS tests exhibited substantially higher power than the Bartlett-type tests for detecting both types of hypotheses (trend and one-point).
 For example, in the case of the one-point-alternative and $\delta\in [0,1.4]$, the tests based on the ATS had about twice as much power than $\varphi_{B-S}$ and $\varphi_{B-P}$ (for which the additional condition is not violated). For the trend-alternative, this power advantage was less pronounced but still clearly visible.  More power simulations on other hypotheses and distributions can be found in the supplementary material.

\section{Review of the required computation time}\label{Time}

Besides power and true type-I-error, the computation time is an important criterion when selecting a proper test.
To take account of this, we performed a small simulation study to compare the computation time for hypothesis $A$ ($\mathcal{H}_0^{\vv}:\vV_{1}=\vV_{2}$) and $B$ ($\mathcal{H}_0^{\vv}:\vV_{1 11}=...=\vV_{1 dd}$).
For each hypothesis and quadratic form both bootstrap techniques were used for 4 different distributions (based on $t_9$-distribution, Normal-distribution, Skew Normal-distribution and Gamma-distribution)  and 2 covariance matrices ($(\vV_1)_{i,j}=0.6^{|i-j|}$ and $\vV_2=\vI_5+\vJ_5$). The average times of 100 such simulation runs are compared. In each run for numerical stability, all of the eight random vectors were considered, and the time is averaged. \\

 For each test 1.000 bootstrap runs were performed with $n_1=125$ observations resp. $\vn=(150,100)$ observations in various dimensions. For the Monte-Carlo-test again, 10.000 simulation steps were used.
  The computations were run by means of the \textsc{R}-computing environment version 3.6.1 \cite{R}  on an Intel Xeon E5430 quad-core CPU
 running at 2.66 GHz using 16 GB DDR2 memory on a Debian GNU Linux 7.8, and the required time in minutes is displayed in \Cref{tab:Zeit1}.
 
 \begin{table}[ht]
\centering
\begin{small}
\begin{tabular}{|l||r|r|r|r|l|r|r|r|r|}\hline
&\multicolumn{4}{|c|}{$A)$}&& \multicolumn{4}{|c|}{$B)$}\\\cline{1-5}\cline{7-10}
  d   & 2&5&10&20  &      & 2&5&10&20   
       \\\cline{1-5}\cline{7-10}
ATS-Para & 0.757 & 5.401 & 27.928 & 222.443 && 0.745 & 5.332 & 27.227 & 195.237  \\\cline{1-5}\cline{7-10}

  ATS-Wild &0.451 & 0.612 & 10.831 & 114.175 && 0.455 & 0.599 & 10.034 & 87.011 \\ \cline{1-5}\cline{7-10}

  ATS &0.086 & 0.195 & 0.387 & 1.231 && 0.056 & 0.154 & 0.276 & 0.698    \\ \cline{1-5}\cline{7-10}

  WTS-Para &0.869 & 6.625 & 44.295 & 462.852 && 0.839 & 6.157 & 34.851 & 275.909 \\ \cline{1-5}\cline{7-10}

  WTS-Wild  &0.559 & 0.916 & 27.252 & 355.679 && 0.545 & 0.794 & 17.681 & 170.063 \\ \cline{1-5}\cline{7-10}

   WTS-$\chi^2$ & 0.003 & 0.004 & 0.038 & 0.320 && 0.003 & 0.003 & 0.029 & 0.133 \\ 
   \hline
 
\end{tabular}
\end{small}
\caption{ Average computation time in seconds of different test-statistics with different dimensions for  hypotheses $A$ ($\mathcal{H}_0^{\vv}:\vV_{1}=\vV_{2}$) and $B$ ($\mathcal{H}_0^{\vv}:\vV_{1 11}=...=\vV_{1 dd}$).}
  \label{tab:Zeit1}
\end{table}
Apart from the classical WTS, all versions of the WTS needed clearly more time than the appropriate ATS. Together with their poor performance in the simulation study, and the additional assumptions on their validity, this makes the WTS unattractive in comparison. 
Moreover, for both the ATS and the WTS, there was a huge difference in the required computation time between the two bootstrap techniques: For small dimensions, the parametric bootstrap needed about 50 percent more computation time than the wild bootstrap, while for larger dimensions, it needed up to more than 20 times longer.
This is not surprising because the generation of normally distributed random vectors is much more time-consuming than generating random weights. 
Moreover, the ATS with the Monte-Carlo based critical values was much faster than all bootstrap approaches as it does not need the repeated calculation of the estimated covariance matrix of the empirical covariances. Additional results on the computation time for other hypotheses can be found in the supplementary material.\\
{\bf Recommendation:} Together with the simulation results, this makes the ATS with parametric bootstrap favorable in the situation with smaller dimensions ($d\leq 5$) due to its accurate type-I-error control. For larger dimensions ($d\geq 10$), however, we recommend its Monte-Carlo implementation due to the much faster computation time.

\subsection{\textsc{Selection of proper hypothesis matrix $\vC$}}\label{Hmatrix}
As mentioned at the beginning, considering a general $\vzeta\neq \vnull_p$ as well as general, not necessarily idempotent and symmetric matrices $\vC$ for the description of the hypothesis is favorable. Beyond more freedom of choosing proper matrices, the major advantage consists of different computational times. 
Indeed,  
depending on the hypothesis of interest, it is possible to choose matrices $\vC\in \R^{m\times ap}$ with $m$ considerably smaller than $ap$.
{We exemplify this issue for the following hypotheses:
\begin{itemize}
\item[$A$)] {Equal Covariance Matrices:\ } Testing the hypothesis $\mathcal{H}_0^{\vv}:\{\vV_{1}=\vV_{2}\} = \{ \vC(A) \vv = \vnull \}$ is usually described by $\vC(A)=\vP_2\otimes \vI_p$. However, the choice $\widetilde\vC(A)=(1,-1)\otimes\vI_p\in \R^{p\times 2p}$ is computationally more efficient.

\item[$B$)] Equal Diagonal Elements: \ The hypothesis  $\mathcal{H}_0^{\vv}:\{\vV_{1 11}=...=\vV_{1 dd}\} = \{ \vC(B) \vv = \vnull \}$ can, e.g., be described by $\vC(B)=\diag(\vh_d)-\vh_d\cdot \vh_d^\top/d$. In contrast, the equivalent description by $\widetilde \vC(B)=(\veins_{d-1},\vnull_{(d-1)\times (d-1)},-\ve_1,\vnull_{(d-1)\times (d-2)},$ $-\ve_2,...,\vnull_{d-1},\ve_{d-1})\in \R^{(d-1)\times p}$ saves a considerable amount of time. Here, $\ve_j$ denotes the $d-1$ dimensional vector containing $1$ in the j-th component and $0$ elsewhere.
\item[$C$)] Equal traces: \ Testing $\mathcal{H}_0^{\vv}:\{\tr(\vV_{1})=\tr(\vV_{2})\} = \{
\vC(C) \vv = \vnull\}$  is usually described by $\vC(C)= \vP_2\otimes [\vh_d\cdot \vh_d^\top]/d$. An equivalent expression is achieved with the smaller matrix $\widetilde\vC(C)=(1,-1)\otimes \vh_d/d\in \R^{1\times 2p}$.
\item[$D$)] Test for a given trace: \ $\mathcal{H}_0^{\vv}:\{\tr(\vV_{1})=\gamma\}=\{
\vC(D) \vv = \vh_d \cdot \gamma \} $ for a given value $\gamma \in \R$ can either be described by  $\widetilde \vC(D)=\vh_d^\top/d\in \R^{1\times p}$ or $\vC(D)= [\vh_d\cdot \vh_d^\top]/d$, where the first choice has considerably less  rows. \end{itemize}}

For these four examples, we performed a small simulation study to compare the computational efficiency of the smaller matrix $\widetilde\vC$ with respect to the quadratic matrix $\vC$: 
 To get reliable results, the same setting as before was used, and the results are displayed in  \Cref{tab:Zeit2} and \Cref{tab:Zeit3}. Depending on the dimension, statistic, and hypothesis of interest, the time savings ranged from less than $1\%$ to more than $99\%$. In fact, for most methods, the savings increased with increasing dimension. Only for the Monte-Carlo ATS test, some fluctuations were visible.
\begin{table}[ht]
\centering
\begin{tabular}{r|c|c|c|c|c|c}
d&\hspace{-0.1cm}$A)$ ATS-Para&\hspace{-0.1cm}$A)$ ATS&\hspace{-0.1cm}$A)$ WTS-Para&\hspace{-0.1cm}$C)$ ATS-Para&\hspace{-0.1cm}$C)$ ATS&\hspace{-0.1cm}$C)$ WTS-Para\\
   \hline
2 & 0.9842 & 0.6516 & 0.9660 & 0.9780 & 0.4299 & 0.9583 \\ \hline
  5 & 0.9872 & 0.7904 & 0.9294 & 0.9713 & 0.2050 & 0.9257 \\ \hline
  10 & 0.9749 & 0.7130 & 0.7868 & 0.9553 & 0.1129 & 0.6893 \\ \hline
 20 & 0.8777 & 0.5669 & 0.5961 & 0.8020 & 0.0966 & 0.4300 \\ 
\end{tabular}
\caption{Computation time for non quadratic hypothesis matrices $\widetilde \vC$ relative to projection matrices $\vC$. Different test-statistics, hypotheses, and dimensions are considered. }
  \label{tab:Zeit2}
\end{table}
\begin{table}[ht]
\centering
\begin{tabular}{r|c|c|c|c|c|c}
d&\hspace{-0.1cm}$B)$ ATS-Para&\hspace{-0.1cm}$B)$ ATS&\hspace{-0.1cm}$B)$ WTS-Para&\hspace{-0.1cm}$D)$ ATS-Para&\hspace{-0.1cm}$D)$ ATS&\hspace{-0.1cm}$D)$ WTS-Para\\
   \hline
2 & 0.8484 & 0.5650 & 0.8493 & 0.8523 & 0.5668 & 0.8465 \\ \hline
 5 & 0.1623 & 0.3811 & 0.1626 & 0.1392 & 0.1888 & 0.1462 \\ \hline
  10 & 0.1014 & 0.4316 & 0.0917 & 0.0336 & 0.1140 & 0.0316 \\ \hline
  20 & 0.0439 & 0.2719 & 0.0372 & 0.0052 & 0.0540 & 0.0040 \\ 
\end{tabular}
\caption{Computation time for non quadratic hypothesis matrices $\widetilde \vC$ relative to projection matrices $\vC$. Different test-statistics, hypotheses, and dimensions are considered. }
  \label{tab:Zeit3}
\end{table}

Moreover, a clear impact of the number of groups could be seen. The reason for this is that for $D)$, the reduction of the dimension can be implemented before the calculation of covariance matrices or similar steps. The latter steps benefitted considerably from this reduction leading to significantly lower computation time.

The exact time measurements for all four hypotheses and both kinds of matrices can be found in the supplementary material.

\section{\textsc{Illustrative Data Analysis}}\label{Illustrative Data Analysis}

To demonstrate the use of the proposed methods, we have re-analyzed neurological data on cognitive impairments. In \cite{bathke2018} the question was  examined whether EEG- or
SPECT-features were preferable to differentiate between three different diagnoses of impairments - 
subjective cognitive complaints  (SCC), mild cognitive impairment (MCI), and Alzheimer disease (AD).
The corresponding trial was conducted at the University Clinic of Salzburg, Department of Neurology. Here one hundred sixty patients were diagnosed with either AD, MCI, or SCC, based on neuropsychological diagnostics, as well as a neurological examination. This data set has been included in the \textsc{R}-package \textit{manova.rm} by \citet*{manova}.
The following \Cref{tab:EEG1} contains the number of patients divided by sex and diagnosis.

\begin{table}[h]
\centering
\begin{tabular}{l|c|c|c|}
&AD&MCI&SCC\\
\hline
male&12&27&20\\
\hline
female&24&30&47\\
\hline
\end{tabular}
\caption{Number of observations for the different factor level combinations of sex and diagnosis.}
\label{tab:EEG1}
\end{table}

For each patient, $d=6$ different kinds of EEG variables were investigated, which leads to  $p=21$ variance and covariance parameters. 
As the male AD and SCC group only contain $12$ and $20$ observations, respectively, an application of the WTS would not be possible.

In \cite{bathke2018}, the authors descriptively checked the empirical covariances matrices to judge that the assumption of equal covariance matrices between the different groups is unlikely. However, this presumption has not been inferred statistically. To close this gap, we first test the null hypothesis of equal covariance matrices between the six different groups using the newly proposed methods. Applying the ATS with parametric resp. wild bootstrap led to p-values of $0.0275$ and $0.0008$.\\
In comparison, the Bartlett-S test of \cite{zhang1993} led to a $p$-value of $0.3484$, potentially reflecting its bad power observed in Section~\ref{Simulations} and also by the authors. Moreover, their Bartlett-P test for the smaller null hypothesis (additionally postulating equality of vectorized moments) shows a small $p$-value of $0.00019998$.

As a next step, we take the underlying factorial structure of the data into account and test, for illustrational purposes, the following hypotheses:

\begin{itemize}
\item[a)] Homogeneity of covariance matrices between different diagnoses,\vspace*{-0.15cm}
\item[b)] Homogeneity of covariance matrices between different sexes,\vspace*{-0.15cm}
\item[c)] Equality of total variance between different diagnosis groups,\vspace*{-0.15cm}
\item[d)] Equality of total variance between different sexes.
\end{itemize}

For the first two hypotheses, we calculated the ATS with wild and parametric bootstrap as well as Bartlett's test statistic with separate and pooled bootstrap. Considering the trace hypotheses, just the first two tests are applicable, and in all cases, the one-sided tests are used based on 10.000 bootstrap runs. The results are presented in \Cref{tab:EEG2} and \Cref{tab:EEG3}.

\begin{table}[h]
\centering
\begin{tabular}{|llc|c|c|c|c|}
\hline
&&& ATS-Para&  ATS-Wild&Bartlett-S&Bartlett-P\\
&&&p-value&p-value&p-value&p-value\\
\hline
$\mathcal{H}_0^a:$&male&AD vs. MCI&0.1000 & \textbf{0.0282} &0.1742 & \textbf{0.0184} \\ \hline
$\mathcal{H}_0^a:$&male&AD vs. SCC& \textbf{<0.0001} & \textbf{<0.0001} & 0.0545 & 0.0634 \\ \hline
$\mathcal{H}_0^a:$&male&MCI vs. SCC & 0.8767 & 0.9801 & 0.1383 & \textbf{0.0078}\\ \hline
$\mathcal{H}_0^a:$&female&AD vs. MCI&  0.0613 & 0.0559 & 0.1050 & 0.1480 \\ \hline
$\mathcal{H}_0^a:$&female&AD vs. SCC& \textbf{0.0128} & \textbf{0.0095} &  \textbf{0.0138} & \textbf{0.0183}\\ \hline
$\mathcal{H}_0^a:$&female&MCI vs. SCC & 0.5656 & 0.6004 & 0.8964 & 0.8988\\ \hline
$\mathcal{H}_0^b:$&AD& male vs. female & 0.1008 & \textbf{0.0279} &0.2479 & 0.0542\\ \hline
$\mathcal{H}_0^b:$&MCI& male vs. female   & 0.2455 & 0.2417 & 0.3695 & 0.4003  \\ \hline
$\mathcal{H}_0^b:$&SCC& male vs. female   & 0.2066 & 0.1914 & 0.2656 & 0.1648\\ \hline

\hline
\end{tabular}
\caption{P-values of ATS with wild resp. parametric bootstrap and Bartlett's test statistic with separate resp. pooled bootstrap for testing equality of covariance matrices. }\label{tab:EEG2}
\end{table}

\begin{table}[h]
\centering
\begin{tabular}{|llc|c|c|}
\hline
&&& ATS-Para&  ATS-Wild\\
&&&p-value&p-value\\
\hline
$\mathcal{H}_0^c:$&male&AD vs. MCI& 0.0733 & 0.0635 \\   \hline
$\mathcal{H}_0^c:$&male&AD vs. SCC& \textbf{<0.0001} & \textbf{<0.0001}\\ \hline
$\mathcal{H}_0^c:$&male&MCI vs. SCC &  0.6146 & 0.6297 \\ \hline
$\mathcal{H}_0^c:$&female&AD vs. MCI&  \textbf{0.0074} & \textbf{0.0091} \\ \hline
$\mathcal{H}_0^c:$&female&AD vs. SCC& \textbf{0.0006} & \textbf{0.0012} \\ \hline
$\mathcal{H}_0^c:$&female&MCI vs. SCC & 0.3687 & 0.3811 \\ \hline
$\mathcal{H}_0^d:$&AD& male vs. female & 0.0881 & 0.0834 \\ \hline
$\mathcal{H}_0^d:$&MCI& male vs. female   & 0.1582 & 0.1592 \\ \hline
$\mathcal{H}_0^d:$&SCC& male vs. female   & 0.3423 & 0.3744  \\ \hline

\end{tabular}
\caption{P-values of ATS with wild resp. parametric bootstrap  for testing equality of traces containing covariance matrices. }\label{tab:EEG3}
\end{table}

It is noticeable that both tests based on the ATS clearly reject the null hypothesis of equal covariances for AD and  SCC for both sexes at level $5\%$, while the p-values of both Bartlett's tests are not significant.
An explanation for this combination with fewer samples may be given by the good small sample performance of the ATS observed in  \Cref{Simulations} and the quite low power of Bartlett's test statistic, which was already mentioned in \cite{zhang1993}. Moreover, the only cases where both Bartlett's test-statistics have smaller p- values are for the combination with the largest sample sizes. Unfortunately, the separate bootstrap has again really low power, while it is questionable whether the additional condition for pooled bootstrap is fulfilled. For the user, this condition is almost as hard to check as equality of covariance. This could lead to the almost circular situation where another test would be necessary to allow for the pooled bootstrap approach for testing homogeneity of covariances.\\

The null hypothesis of equal total covariance resp. equal traces could be rejected significantly (at level 5\%) by both bootstrap tests in three cases. 
Perhaps surprising at first is that the null hypothesis of equal covariance matrices between the female AD and MCI groups could not be rejected, but the joint univariate null hypothesis of equal traces could now be rejected at level $5\%$. 

Although the hypothesis of equal covariance matrices couldn't be rejected in each case, it shows that sex and diagnosis are likely to have an effect on the covariance matrix. This illustrative analysis underpins that the approach of \cite{bathke2018}, which can deal with covariance heterogeneity, was very reasonable.

\section{\textsc{Conclusion \& Outlook}}\label{Conclusion & Outlook}

In the present paper, we have introduced and evaluated a unified approach to testing a variety of quite general null hypotheses formulated in terms of covariance matrices. 
The proposed method is valid under a comparatively small number of requirements that are verifiable in practice. 
Previously existing procedures for the situation addressed here had suffered from low power to detect alternatives, were limited to only a few specific null hypotheses, or needed various requirements in particular regarding the data generating distribution. \\
Under weak conditions, we have proved the asymptotic normality of the difference between the vectorized covariance matrices and their corresponding vectorized empirical versions. 
We considered two-test statistics, which are based upon the vectorized 
empirical covariance matrix and an estimator of its own covariance:  a 
Wald-type-statistic (WTS) as well as an ANOVA-type-statistic (ATS).
These exhibit the usual advantages and disadvantages that are already well-known from the literature on mean-based inference. 
In order to take care of some of these difficulties, namely the critical value for the ATS being unknown and the WTS requiring a rather large sample size, two kinds of bootstrap were used. 
On this occasion, specific adaptions were needed to take account of the special situation where inference is not on the expectation vectors but on the covariance matrices. \\
To investigate the properties of the newly constructed tests, an extensive simulation study was done. 
For this purpose, several different hypotheses were considered, and the type-I-error control, as well as the power to detect deviations from the null hypothesis, were compared to existing test procedures.
The ATS showed a quite accurate error control in each of the hypotheses, in particular in comparison with competing procedures.  Note that for most hypotheses, no appropriate competing test is available. 
The simulated power of the proposed tests was fine, even for moderately small sample sizes ($n_1=30,n_2=20$). 
This is a major advantage when comparing with existing procedures for testing homogeneity of covariances, even considering that they usually require further assumptions .\\\\
In future research, we would like to investigate in more detail the large number of possible null hypotheses that are included in our model as special cases. For example, tests for given covariance structures (such as compound symmetry or autoregressive) with unknown parameters are of great interest.
Moreover, our results allow for a variety of new tests for hypotheses that can be derived from our model, for example, testing the equality of determinants of covariances matrices.
The model and the assumption of finite fourth moments exclude some distributions like, for example, heavy-tailed distributions. For such distributions, probably a similar approach can be developed using scatter matrices.\\
 Finally, it is still unclear whether our approach can be extended to high-dimensional settings. There already exist some inspiring solutions, see for example, \cite{chi2012global}, \cite{li2012}, \cite{li2014}), and \cite{cai2013}.
However, they are only constructed for special situations and do not allow the same flexibility as our approach. Due to different technical approaches, this task remains future research.
Furthermore, we are planning to investigate extensions of our work by combining it with results on high-dimensional covariance matrix estimators, as considered in \citet*{cai2016}.

\section{\textsc{Acknowledgment}}
Paavo Sattler and Markus Pauly would like to thank the German Research Foundation for the support received within project PA 2409/4-1. Moreover, Arne Bathke expresses his thanks to the Austrian Science Fund (FWF) for the funding received through project I 2697-N31.

\section{Appendix}
\subsection{The Model }\label{The Model}
\lhead[\footnotesize\thepage\fancyplain{}\leftmark]{}\rhead[]{\fancyplain{}\rightmark\footnotesize\thepage}
{The considered semiparametric model can be shortly defined through $\vX_{ik}=\mu_i+\vepsilon_{ik}$, while $\vep_{ik}$ are i.i.d. d-dimensional random vectors with $\E(\vep_{ik})=\vnull_d$ and $\Var(\vep_{ik})=\vV_i\geq 0$. For $\vep_{ik}$, which is called the non-parametric part, it is allowed that each component are from a completly different distribution. We additionaly assume finite fourth moments for all components through  $\E(||\vep_{ik}||^4)<\infty$. The number of groups $a$ can persist from multiple crossed factors, where for one, two and three factors the model is given through:
\[\begin{array}{ll}
\vX_{ik}&=\vmu_i+\vep_{ik}\\
&=\vmu+\valpha_i+\vbeta_j+(\valpha\vbeta)_{ij}+\vep_{ijk}\\
&=\vmu+\valpha_i+\vbeta_j+\vgamma_\ell+(\valpha\vbeta)_{ij}+(\vbeta \vgamma)_{j\ell}+(\valpha \vgamma)_{i\ell}+(\valpha \vbeta \vgamma)_{ij\ell}+\vep_{ij\ell k},\\
\end{array}\]
 with $i=1,...,a$ , $j=1,...,J$ , $\ell=1,....,L$ and $k=1,...,K$. In this case it holds $a=J\cdot L\cdot K$, and each group represents one combination of these three factors.\\
Like this, the dimension can also consist of multiple factors, for example, in repeated measure design, where often the factors time and treatment are crossed. Then the dimension is divided into smaller parts, one for each factor-combination. }

\subsection{Proofs }
The asymptotic distribution,  discussed in \Cref{Theorem1} is well known (for example, from \cite{browne}), but based on the importance of the techniques presented in this paper, we will prove it shortly. Moreover, this allows getting the idea of our bootstrap approaches later on.
\begin{proof}[Proof of \Cref{Theorem1}]

First we consider the difference between the vector $\vv_i$ and its estimated version $\widehat \vv_i$, multiplied with $\sqrt{N}$
\[\begin{array}{c} {\sqrt{N}}(\widehat \vv_i  -\vv_i)\\= {\sqrt{N}}\vech\left(\frac 1 {n_i-1} \sum\limits_{k=1}^{n_i}\left[\vep_{i k}\vep_{i k}^\top-{\vV}_i\right]+ \frac 1 {n_i-1}{\vV}_{i}-\frac {1}{n_i-1}(\sqrt{n_i} \ \overline{\vep}_{i \cdot})(\sqrt{n_i} \ \overline{\vep}_{i \cdot})^\top \right).\end{array}\]

Due to  Slutzky and the multivariate Central limit theorem, the second and third term tends to zero in probability. Thus,  it is sufficient to consider the first term.  But this converges to $\mathcal{N}_{ d}(\vnull_d,\kappa_i^{-1}\vSigma_i)$ in distribution again by the multivariate central limit theorem, which gives us the result due to independence of the groups. 
\end{proof}
This convergence would also follow from \cite{zhang1993}, but the bootstrap approach is based on this proof, so it is helpful to outline it again. To use this result, a consistent estimator for the covariance matrix $\vSigma$ is needed.

\begin{proof}[Consistency of $\widehat \vSigma$]
 Because $\vech(\vep_{ik}\vep_{ik}^\top)$ are i.i.d. vectors, we know that \[\widetilde \vSigma_i=\frac{\sum\limits_{k=1}^{n_i}\left[\vech(\vep_{ik}\vep_{ik}^\top)-\sum\limits_{\ell=1}^{n_i}\frac{\vech(\vep_{i\ell}\vep_{i\ell}^\top)}{n_i}\right]\left[\vech(\vep_{ik}\vep_{ik}^\top)-\sum\limits_{\ell=1}^{n_i}\frac{\vech(\vep_{i\ell}\vep_{i\ell}^\top)}{n_i}\right]^\top}{n_i-1}\]
 is a consistent estimator for $\vSigma_i$.
However, we can not calculate this estimator, since $\Cov(\vep_i)=\Cov(\vX_i)$, but in general $\Cov(\vech(\vep_i\vep_i^\top))\neq \Cov(\vech(\vX_i\vX_i^\top))$ and $\vmu_i$ is unknown. So we use  centered vectors $\widetilde \vX_{ik}$ to formulate the proper covariance matrix $\widehat \vSigma_i$. These vectors are not independent, so we prove the consistency of our estimator through that $\widehat \vSigma_i - \widetilde \vSigma_i$ converge almost sure to 0.  
 This is done by\\\\
 $\begin{array}{lll}&&\hspace*{-0.15cm}\widehat \vSigma_i - \widetilde \vSigma_i\\
=&&\hspace*{-0.15cm} \frac{4n_i}{n_i-1} \left(\vech(\overline \vX_i \vmu_i^\top) \vech(\overline \vX_i \vmu_i^\top) ^\top-\vech(\overline \vX_i \overline \vX_i^\top) \vech(\overline \vX_i \overline \vX_i ^\top) ^\top\right)\\
 &\hspace*{-0.15cm}+&\hspace*{-0.15cm} \frac{4}{n_i-1}\sum_{k=1}^{n_i} \left[ \vech( \vX_{ik} \overline \vX_i ^\top) \vech( \vX_{ik} \overline \vX_i ^\top) ^\top-\vech( \vX_{ik} \vmu_i^\top) \vech( \vX_{ik} \vmu_i^\top) ^\top\right]\\[1.2ex]
 &\hspace*{-0.15cm}+&\hspace*{-0.15cm}\frac{4}{n_i-1}\sum_{k=1}^{n_i}  \left[\vech(\vX_{ik} \vX_{ik}^\top) \vech(\overline \vX_i \vmu_i^\top) ^\top-\vech(\vX_{ik} \vX_{ik}^\top) \vech(\overline \vX_i \overline \vX_i ^\top) ^\top\right]\\[1.2ex]
 &\hspace*{-0.15cm}+&
\hspace*{-0.15cm}\frac{4}{n_i-1}\sum_{k=1}^{n_i}  \left[\vech(\vX_{ik} \vX_{ik}^\top) \vech( \vX_{ik} \vmu_i^\top)  ^\top-\vech(\vX_{ik} \vX_{ik}^\top) \vech( \vX_{ik} \overline \vX_i ^\top)  ^\top\right] \end{array}$\\
$\begin{array}{lll}

  =&& \hspace*{-0.15cm}\frac{4n_i}{n_i-1} \left(\vech(\overline \vX_i (\vmu_i-\overline \vX_i)^\top) \vech(\overline \vX_i \vmu_i^\top+\overline \vX_i \overline \vX_i ^\top)^\top\right)\\
 &\hspace*{-0.15cm}+&\hspace*{-0.15cm}\frac{4}{n_i-1}\sum_{k=1}^{n_i} \left[ \vech( \vX_{ik} (\overline \vX_i - \vmu_i)^\top) \vech( \vX_{ik} (\overline \vX_i-\vmu_i) ^\top+2\vX_{ik} \vmu_i^\top) ^\top\right]\\
 &\hspace*{-0.15cm}+&\hspace*{-0.15cm}\frac{4}{n_i-1}\sum_{k=1}^{n_i}  \left[\vech(\vX_{ik} \vX_{ik}^\top) \vech(\overline \vX_i (\vmu_i-\overline \vX_i)^\top) ^\top\right]\\&\hspace*{-0.15cm}+&
\hspace*{-0.15cm}\frac{4}{n_i-1}\sum_{k=1}^{n_i}  \left[\vech(\vX_{ik} \vX_{ik}^\top) \vech( \vX_{ik} (\vmu_i-\overline \vX_i)^\top)  ^\top\right].
 \end{array}$\\\\\\
  It is enough  to show that each component of this difference converges almost sure to zero. So with $|\vX|$ denoting the absolute value of each component we get  for arbitrary $h,j\in \{1,...,p\}$ that \\\\
  $\begin{array}{lll}&& \hspace*{-0.15cm}|(\widehat \vSigma_i - \widetilde \vSigma_i)_{h,j}|\\
\leq&& \hspace*{-0.15cm}\frac{4n_i}{n_i-1} \lvert\vech(\overline \vX_i (\vmu_i-\overline \vX_i)^\top)_h\lvert \cdot \lvert \vech(\overline \vX_i \vmu_i^\top+\overline \vX_i \overline \vX_i ^\top)_j\lvert\\
 &\hspace*{-0.15cm}+&\hspace*{-0.15cm}\frac{4}{n_i-1}\sum_{k=1}^{n_i}\lvert  \vech( \vX_{ik} (\overline \vX_i - \vmu_i)^\top)_h\lvert \cdot  \lvert \vech( \vX_{ik} (\overline \vX_i-\vmu_i) ^\top+2\vX_{ik} \vmu_i^\top) _j\lvert\\
 &\hspace*{-0.15cm}+&\hspace*{-0.15cm}\frac{4}{n_i-1}\sum_{k=1}^{n_i}  \lvert \vech(\vX_{ik} \vX_{ik}^\top)_j \lvert \cdot \lvert \vech(\overline \vX_i (\vmu_i-\overline \vX_i)^\top)_h \lvert\\&\hspace*{-0.15cm}+&
\hspace*{-0.15cm}\frac{4}{n_i-1}\sum_{k=1}^{n_i}  \lvert\vech(\vX_{ik} \vX_{ik}^\top)_j\lvert \cdot \lvert \vech( \vX_{ik} (\vmu_i-\overline \vX_i)^\top)_h\lvert
 \end{array}$\\\\
  $\begin{array}{lll}
\leq&&  \hspace*{-0.15cm}\max_{\ell=1,...p} \lvert(\vmu_i)_\ell-(\overline {\vX_i})_\ell\lvert \cdot \frac{4n_i}{n_i-1}\cdot  \vech( \lvert\overline {\vX_i}\lvert \veins^\top)_h \cdot \lvert \vech(\overline \vX_i \vmu^\top+\overline \vX_i \overline \vX_i ^\top)_j\lvert\\[0.5ex]
 &\hspace*{-0.15cm}+&\hspace*{-0.15cm}\left(\max_{\ell=1,...p} \lvert(\vmu_i)_\ell-(\overline {\vX_i})_\ell\lvert\right)^2\cdot \frac{4}{n_i-1} \sum_{k=1}^{n_i}  \vech( \lvert\vX_{ik}\lvert \veins^\top)_h \cdot   \vech(\lvert \vX_{ik}\lvert \veins^\top) _j\\
 &\hspace*{-0.15cm}+&\hspace*{-0.15cm}\max_{\ell=1,...p} \lvert(\vmu_i)_\ell-(\overline {\vX_i})_\ell\lvert\cdot \frac{4}{n_i-1} \sum_{k=1}^{n_i}  \vech( \lvert\vX_{ik}\lvert \veins^\top)_h \cdot  \lvert \vech( 2\vX_{ik} \vmu^\top) _j\lvert\\ 
 
 &\hspace*{-0.15cm}+&\hspace*{-0.15cm}\max_{\ell=1,...p} \lvert(\vmu_i)_\ell-(\overline {\vX_i})_\ell\lvert\cdot \frac{4}{n_i-1}\sum_{k=1}^{n_i}  \lvert \vech(\vX_{ik} \vX_{ik}^\top)_j \lvert \cdot  \vech(\lvert\overline \vX_i\lvert \veins^\top)_h \\&\hspace*{-0.15cm}+&\hspace*{-0.15cm}
\max_{\ell=1,...p} \lvert (\vmu_i)_\ell-(\overline {\vX_i})_\ell\lvert\cdot \frac{4}{n_i-1} \sum_{k=1}^{n_i}  \lvert\vech(\vX_{ik} \vX_{ik}^\top)_j\lvert \cdot  \vech(\lvert \vX_{ik}\lvert \veins^\top)_h.
 \end{array}$\\ \\
 Here we used that the maximum doesn't depend on the index of the sum, so this factor can be pulled out of the $\vech$ and the sum, which are both linear functions.
Because of the strong law of large numbers we know  $(\vmu-\overline \vX_i)\stackrel{a.s.}{\to} \vnull_d$ which means that every component goes to zero almost sure and therefore also the maximum.   \\\\
The general assumption (4) , which ensures that all occurring terms have finite expectation values together  with another application of the SLLN  leads to:
 \[\vech( \lvert\overline {\vX_i}\lvert \veins^\top)_h \cdot \lvert \vech(\overline \vX_i \vmu^\top+\overline \vX_i \overline \vX_i ^\top)_j\lvert\stackrel{a.s.}{\longrightarrow} \vech( \lvert\mu_i\lvert \veins^\top)_h \cdot \lvert \vech(2\mu_i \vmu_i^\top)_j\lvert,\]
 \[\frac{4}{n_i-1} \sum_{k=1}^{n_i}  \vech( \lvert\vX_{ik}\lvert \veins^\top)_h \cdot   \vech(\lvert \vX_{ik}\lvert \veins^\top) _j\stackrel{a.s.}{\longrightarrow}
4 \cdot \E\left(  \vech( \lvert\vX_{i1}\lvert \veins^\top)_h \cdot   \vech(\lvert \vX_{i1}\lvert \veins^\top) _j\right)
 \]\\
and equivalent for the other sums. So we have in all this cases the products goes almost sure to zero and therefore $|(\widehat \vSigma_i - \widetilde \vSigma_i)_{h,j}|\stackrel{a.s.}{\longrightarrow} 0$. Because of the independence of the groups, we also get the result for $\widehat\vSigma$.
\end{proof}

With these results, the asymptotic distribution of the applied test statistics can be prooved.

\begin{proof}[Proof of \Cref{Verteilung}]
All results are known (see, e.g., \cite{BrunnerBathkeKonietschke}), but sometimes only idempotent symmetric hypothesis matrices are considered, so we will repeat them for general matrices $\vC$.
From \Cref{Theorem1} it follows that all these quadratic forms can be written as the sum of a quadratic form with normal distributed random vectors and vectors which converge in distribution to zero. 

Therefore with $\vZ\sim\mathcal{N}_{a\cdot p}\left(\vnull_{a\cdot p},\vI_{ap} \right)$ and $\sqrt{N} (\widehat \vv -\vv)\stackrel{\mathcal {D}}{\longrightarrow}\vSigma^{1/2} \vZ$ we get
\[\begin{array}{ll}\widehat Q_{\vv}&= N\left[  \vC\widehat \vv - \vzeta\right]^\top  \vE(\vC,\widehat{\vSigma})\left[  \vC\widehat \vv - 
\vzeta\right]\\&\stackrel{\mathcal H_0}{=}{N}\cdot (\widehat \vv -\vv)^\top\vC^\top  \vE(\vC,\widehat{\vSigma}) \vC(\widehat \vv -\vv)
\\&\stackrel{\mathcal {D}}{\rightarrow} \left(\vSigma^{1/2}\vZ\right)^\top \vC^\top  \vE(\vC,\widehat{\vSigma})\vC \left(\vSigma^{1/2} \vZ\right)
\\&=\vZ^\top \vSigma^{1/2}\vC^\top   \vE(\vC,\widehat{\vSigma}) \vC \vSigma^{1/2}  \vZ\\&\stackrel{\mathcal {D}}{=} \sum_{\ell=1}^{a\cdot p} \lambda_\ell B_\ell,\end{array}\]
with $\lambda_\ell, \ell=1,...,ap$ eigenvalues of$(\vSigma^{1/2}\vC^\top\vE(\vC,\vSigma)\vC \vSigma^{1/2})$ and $B_\ell\stackrel{i.i.d.}{\sim} \chi_1^2$.
{Note, that we have used that $(\vSigma^{1/2}\vC^\top\vE(\vC,\vSigma)\vC \vSigma^{1/2})$ is symmetric and therefore has a spectral representation. The rest of the proof follows from the fact that the multivariate standard normal distribution is invariant under orthogonal transformations, the consistency of $\vE(\vC,\widehat \vSigma)$ for $\vE(\vC, \vSigma)$ and the continuous mapping theorem.}\end{proof}

\begin{proof}[Proof of \Cref{PBTheorem1}]
It is sufficient to prove the part for the single groups because the second part is just the combination of all groups.\\
This result follows from a part-wise application (given the data) of the multivariate Lindeberg-Feller-Theorem. So it remains to show that all conditions are fulfilled, for which we use  the fact that $\vY^*$ under $\vX$ is $p$ dimensional normal distributed with  expectation $\vnull_p$ and variance  $\widehat{\vSigma}_{i}$:\\

$\begin{array}{lll}1.)\hspace*{-0.3cm}&\textcolor{white}{=} &\sum\limits_{k=1}^{n_i}\E\left(\frac {\sqrt{N}} {n_i}\vY_{i k}^*\Big\lvert \vX\right)=\sum\limits_{k=1}^{n_i}\frac {\sqrt{N}} {n_i}\cdot \E\left(\vY_{i k}^*\Big\lvert \vX\right)=0\end{array}$\\\\
$\begin{array}{lll}2.)\hspace*{-0.3cm}&\textcolor{white}{=} &\sum\limits_{k=1}^{n_i}\Cov\left(\frac {\sqrt{N}} {n_i}\vY_{i k}^*\Big\lvert \vX\right)= \sum\limits_{k=1}^{n_i}\frac {{N}} {n_i^2}  \widehat\vSigma_{i}\stackrel{\mathcal P}{\to} \frac 1 {\kappa_i} \vSigma_{i}\end{array}$\\

$\begin{array}{lll}3).\hspace*{-0.3cm}&& \lim\limits_{N\to \infty}\sum\limits_{k=1}^{n_i} \E\left( \Big \lvert \Big \lvert \frac {\sqrt{N}} {n_i}\vY_{i k}^*  \Big \lvert \Big \lvert^2\cdot \ind_{ \big \lvert \big \lvert \frac {\sqrt{N}} {n_i}\vY_{i k}^*  \big \lvert \big \lvert>\delta}\ \Big\lvert \vX\right)\\[1.4ex]

&=& \lim\limits_{N\to \infty}  \frac {{N}} {n_i^2}\sum\limits_{k=1}^{n_i} \E\left( \big  \lvert \big  \lvert \vY_{i 1}^*  \big  \lvert \big  \lvert^2\cdot \ind_{ \lvert \lvert \vY_{i 1}^*   \lvert  \lvert>\delta \frac {n_i}{\sqrt{N}} }\ \Big\lvert \vX\right)\\[1.6ex]
&=& \frac 1 {\kappa_i}  \cdot \lim\limits_{N\to \infty}  \E\left( \big  \lvert \big  \lvert \vY_{i 1}^*\big \lvert \big \lvert^2\cdot \ind_{  \lvert \lvert \vY_{i 1}^* \lvert \lvert>\delta\frac {n_i}{\sqrt{N}} }\ \Big\lvert \vX\right)\\[1.6ex]
&\leq &\frac 1 {\kappa_i}  \cdot \lim\limits_{N\to \infty} \sqrt{ \E\left(  \lvert   \lvert \vY_{i 1}^* \lvert \lvert^2 \ \lvert \vX\right)}\cdot \sqrt{\E\left(\ind_{  \lvert  \lvert \vY_{i 1}^*  \lvert  \lvert>\delta\frac {n_i}{\sqrt{N}} }\ \Big\lvert \vX\right)}=0 

\end{array}$\\\\\\
Here we used the Cauchy-Bunjakowski-Schwarz-Inequality and that we know $\E\left( \big  \lvert \big  \lvert \vY_{i 1}^*\big \lvert \big \lvert^2\ \big\lvert \vX\right)$. Moreover  because of  the condition  $ {n_i}/ N \to \kappa_i$ and as a consequence {$\delta\cdot  {n_i}/{\sqrt{N}}\to \infty $} it holds $P\left( \big \lvert \big \lvert \vY_{i 1}^* \big \lvert \big \lvert>\delta \cdot {n_i}/\sqrt{N} \right)\to 0$, which leads to the result.\\

Therefore, given the data $\vX$
it follows that ${\sqrt{N}}\cdot \overline \vY_{i }^*$  converges in distribution to $\mathcal{N}_{p}\left(\vnull_p, 1/ {\kappa_i}\cdot  \vSigma_i\right)$ and due to indepence also
$ {\sqrt{N}}\cdot \overline \vY^*$ converges in distribution to $\mathcal{N}_{a\cdot p}\left(\vnull_{a\cdot p},\bigoplus_{i=1}^a 1/ {\kappa_i}\cdot  \vSigma_i\right)$.
{As the empirical covariance matrix of the bootstrap sample is also consistent, with {$\widehat \vSigma_i^*\stackrel{\mathcal{P}}{\to}\widehat \vSigma_i$} the result follows from  Consistency of $\widehat \vSigma_i$ and the triangle inequality. Moreover, $\widehat \vSigma^*\stackrel{\mathcal{P}}{\to} \vSigma$ follows by continuous mapping theorem.}\end{proof}

\begin{proof}[Proof of \Cref{WBTheorem}]
 Again we have to show the conditions of the Lindeberg-Feller theorem part-wise, given the data  $\vX=(\vX_{11}^\top, \dots, \vX_{an_a}^\top)^\top$ :\\\\
$\begin{array}{lll}1.)\hspace*{-0.3cm}& \textcolor{white}{=}&\sum\limits_{k=1}^{n_i}\E\left(\frac {\sqrt{N}} {n_i}\vY_{i k}^\star\Big\lvert \vX\right)=\sum\limits_{k=1}^{n_i}\frac {\sqrt{N}} {n_i}\E(W_{ik})\cdot \left[\vech(\widetilde \vX_{ik}\widetilde \vX_{ik}^\top )- \sum\limits_{i=1}^{n_i}\frac  {\vech(\widetilde \vX_{ik}\widetilde \vX_{ik}^\top )}{n_i}\right] =0\end{array}$\\\\
$\begin{array}{lll}2.)\hspace*{-0.3cm}& &\sum\limits_{k=1}^{n_i}\Cov\left(\frac {\sqrt{N}} {n_i}\vY_{i k}^\star\Big\lvert \vX\right)= \frac {{N}} {n_i^2} \E\left(W_{i1}^2\right)\cdot(n_i-1)\cdot \widehat \vSigma_i\\ [1.6ex]
& =& \frac{n_i-1}{n_i}\frac {{N}} {n_i} \widehat\vSigma_{i}\stackrel{\mathcal P}{\to} \frac 1 {\kappa_i} \vSigma_{i}\end{array}$\\\\\\
For the last part we use  that given the data $\big  \lvert \big  \lvert \vY_{i 1}^\star  \big  \lvert \big  \lvert^2\cdot \ind_{  \lvert  \lvert \vY_{i 1}^\star  \lvert  \lvert>\delta \frac {n_i}{\sqrt{N}} }\leq \big  \lvert \big  \lvert \vY_{i 1}^\star  \big  \lvert \big  \lvert^2$ has a finite expectation value. Moreover  Lebesgue's dominated convergence theorem with $ {n_i}/{\sqrt{N}}\to \infty$ and $P\left( \big \lvert \big \lvert \vY_{i 1}^\star \big \lvert \big \lvert>\delta \cdot {n_i}/\sqrt{N} \right)\to 0,$ leads to the result.\\

$\begin{array}{lcl}3).&\quad&\lim\limits_{n_i\to \infty}\sum\limits_{k=1}^{n_i} \E\left( \Big \lvert \Big \lvert \frac {\sqrt{N}} {n_i}\vY_{i k}^\star  \Big \lvert \Big \lvert^2\cdot \ind_{ \big \lvert \big \lvert \frac {\sqrt{N}} {n_i}\vY_{i k}^\star  \big \lvert \big \lvert>\delta}\ \Big\lvert \vX\right)\\[1.3ex]

&=& \lim\limits_{N\to \infty}  \frac {{N}} {(n_i)^2}\sum\limits_{k=1}^{n_i} \E\left( \big \lvert \big  \lvert \vY_{i 1}^\star \big  \lvert \big  \lvert^2\cdot \ind_{  \lvert  \lvert \vY_{i 1}^\star   \lvert  \lvert>\delta \frac {n_i}{\sqrt{N}} }\ \Big \lvert \vX\right)\\[1.3ex]
&=& \frac 1 {\kappa_i} \cdot \lim\limits_{N\to \infty}  \E\left( \big  \lvert \big  \lvert \vY_{i 1}^\star  \big  \lvert \big  \lvert^2\cdot \ind_{ \lvert  \lvert \vY_{i 1}^\star  \lvert  \lvert>\delta \frac {n_i}{\sqrt{N}} }\ \Big \lvert \vX\right)\\[1.3ex]
&=& \frac 1 {\kappa_i}  \cdot    \E\left(\lim\limits_{N\to \infty}   \big  \lvert \big  \lvert \vY_{i 1}^\star  \big  \lvert \big  \lvert^2\cdot \ind_{  \lvert  \lvert \vY_{i 1}^\star  \lvert  \lvert>\delta \frac {n_i}{\sqrt{N}} }\ \Big \lvert \vX\right)
=0    
\end{array}$\\\\
Hence, given the data we have convergence in distribution of
${\sqrt{N}}\cdot \overline \vY_{i }^\star$  and $ {\sqrt{N}}\cdot \overline \vY^\star$ to $\mathcal{N}_{p}\left(\vnull_p, 1/ {\kappa_i}\cdot  \vSigma_i\right)$ resp. $\mathcal{N}_{a\cdot p}\left(\vnull_{a\cdot p},\bigoplus_{i=1}^a 1/ {\kappa_i}\cdot  \vSigma_i\right)$.\\
The consistency of the covariance estimator is proven 
analogous to the parametric bootstrap. \end{proof}

 \begin{proof}
[Proof of \Cref{KorParametric} and \Cref{KorWild} ]

{As in \Cref{Verteilung} it holds that
\[ N\left[  \vC\widehat \vv - \vzeta\right]^\top  \vE(\vC,\vSigma)\left[  \vC\widehat \vv - 
\vzeta\right]\stackrel{\mathcal {D}}{\To} \sum_{\ell=1}^{a\cdot p} \lambda_\ell B_\ell,\]
where $\lambda_\ell, \ell=1,...,ap$ are the eigenvalues of $(\vSigma^{1/2}\vC^\top\vE(\vC,\vSigma)\vC \vSigma^{1/2})$ and $B_\ell\stackrel{i.i.d.}{\sim} \chi_1^2$.
Moreover, similar to \Cref{WBTheorem} it follows that given the data, \[N\left[ \vC\overline \vY^*  \right]^\top \vE(\vC,\widehat \vSigma^*)\left[ \vC\ \overline \vY^*\right]\stackrel{\mathcal {D}}{\to} \sum_{\ell=1}^{a\cdot p} \lambda_\ell B_\ell\] and
\[N\left[ \vC\overline \vY^\star  \right]^\top \vE(\vC,\widehat \vSigma^\star)\left[ \vC\ \overline \vY^\star\right]\stackrel{\mathcal {D}}{\to} \sum_{\ell=1}^{a\cdot p} \lambda_\ell B_\ell,\]
because $\widehat{\vSigma}^*$ and $\widehat{\vSigma}^\star$ are consistent estimators for ${\vSigma}$. }\end{proof}
The result especially  allows the application of the parametric bootstrap version of the MATS given by
\begin{equation}\label{MATSP} MATS^*:=N \left[\vC\overline \vY^*\right]^\top \left(\vC \widehat\vSigma_0^*\vC^\top\right)^+ \left[\vC\overline \vY^*\right]\end{equation}
and the  wild bootstrap version given by
\begin{equation}\label{MATSW} MATS^\star:=N \left[\vC\overline \vY^\star\right]^\top \left(\vC \widehat \vSigma_0^\star\vC^\top\right)^+ \left[\vC\overline \vY^\star\right].\end{equation}

\subsection{\textsc{Further Simulations} }
{ In this section, we expand the simulations from \Cref{Simulations}, for example, through more null hypotheses and bootstrap versions of the MATS statistic defined in \eqref{MATSP} and \eqref{MATSW}. To investigate the influence of the covariance matrix, for the distributional setting an additional covariance matrix $\vV_2$ is used, which is a compound symmetry matrix given by $\vV_2:=\vI_5+\vJ_5$. The same distributions as in \Cref{Simulations} were used for the error term, but we also simulated one more, which is based on a standardized centered t-distribution with 9 degrees of freedom. }

{Testing for the equality of covariances is an important hypothesis that usually becomes more demanding for an increasing number of groups. Therefore, for all the random vectors, we investigated an additional scenario:}
\begin{itemize}
{\item[$E$)]$a=3$\quad $\mathcal{H}_0^{\vv}:\vV_{1}=\vV_{2}=\vV_3$,}

\end{itemize} {where also scenario $E)$ can be formulated with an idempotent symmetric matrix $\vC(E)=\vP_3\otimes \vI_{15}$. For scenario $A)$ and $C)$ we considered $n_1= 0.6\cdot N$ and $n_2=0.4\cdot N$ with {$\vN=(50,100,250,500)$} and for  $B)$ {$\vn_1=(25,50,125,250)$}.
In case of the three groups we considered $n_1:=0.4\cdot N$, $n_2:=0.25\cdot N$ and $n_3:=0.35\cdot N$ for $N$ from 80 up to 800. This choice makes the sample sizes similar to the situation with two groups and therefore increases the comparability.}

{
We should keep in mind that in this case, p is $15$, which makes some of these sample sizes small in relation to the dimension.
The $WTS$ resp. $MATS$ are part of our simulation,  although in practice it is quite difficult or even impossible to check the necessary conditions $\vSigma>0$ resp. $\vSigma_0 >0$.}

{
Again it could be seen in all tables that the wild bootstrap lead to more liberal results and the parametric bootstrap had less liberal or even conservative test results. This hold for all our quadratic forms, the ATS, the WTS, and the MATS. Overall hypotheses and settings, the MATS-test-statistic seems to perform between the ATS and the WTS but was still preferable over the Bartlett test-statistics in scenario A).

\renewcommand{\baselinestretch}{1}
\begin{table}[htbp]
\centering
\begin{tabular}{|l||c|c|c|c||c||c|c|c|c||}\hline
    &\multicolumn{4}{|c||}{$t_9$}&&\multicolumn{4}{|c||}{Normal}\\\cline{1-5} \cline{7-10}
 \hspace*{.1cm}N&   50&100&250&500&&   50&100&250&500
       \\\cline{1-5} \cline{7-10}
ATS-Para & {\bf.0494}& {\bf.0525} & {\bf.0496} & {\bf.0504 }&& .0579 & .0540 & {\bf.0518} & {\bf.0515 } \\ \cline{1-5} \cline{7-10}
  ATS-Wild&.0792 & .0698 & .0580 & .0553 && .0797 & .0672 & .0558 & .0533 \\ \cline{1-5} \cline{7-10}
  ATS &.0552 & .0537 & {\bf.0498} & {\bf.0498} && .0634 & .0562 & {\bf.0520} & {\bf.0510}  \\ \cline{1-5} \cline{7-10}
  WTS-Para& .0627 & .0638 & .0596 & {\bf.0547} && .0659 & .0661 & .0623 & .0566  \\ \cline{1-5} \cline{7-10}
  WTS-Wild& .0980 & .0895 & .0726 & .0643 && .0961 & .0852 & .0706 & .0612  \\ \cline{1-5} \cline{7-10}
  WTS-$\chi_{15}^2$  &.4965 & .2168 & .1002 & .0738 && .5000 & .2161 & .1054 & .0757  \\ \cline{1-5} \cline{7-10}
  MATS-Para & .0594 & .0598 & .0546 & {\bf.0525} && .0649 & .0596 & .0538 & .0534  \\ \cline{1-5} \cline{7-10}
  MATS-Wild&.0838 & .0724 & .0604 & .0554 && .0853 & .0694 & .0576 & .0553 \\ \cline{1-5} \cline{7-10}
  Bartlett-S & .0168 & {\bf.0492} & .0577 & {\bf.0524} && .0111 & .0371 &{\bf .0478} & {\bf.0485}\\ \cline{1-5} \cline{7-10}
  Bartlett-P & .0233 & .0392 & .0464 & .0465 && .0199 & .0360 & .0452 & .0467 \\ \cline{1-5} \cline{7-10}
  Box's M-$\chi_{15}^2$& .1308 & .1337 & .1361 & .1401 && .0638 & .0575 & {\bf.0521} & {\bf.0496} \\ \cline{1-5} \cline{7-10}
  Box's M-F& .1238 & .01322 & .1358 & .1400 && .0609 & .0567 & {\bf.0520} & {\bf.0496} \\ \cline{1-5} \cline{7-10}
   \hline\multicolumn{10}{|c|}{}\\\hline
   &\multicolumn{4}{|c||}{Skew Normal}&&\multicolumn{4}{|c||}{Gamma}\\\cline{1-5} \cline{7-10}
   \hspace*{.1cm}N&   50&100&250&500&&   50&100&250&500
       \\\cline{1-5} \cline{7-10}
ATS-Para & .0589 & .0538 &{\bf .0528} & {\bf.0488} && {\bf.0485} & .0439 & .0439 & .0464\\ \cline{1-5} \cline{7-10}
  ATS-Wild& .0915 & .0708 & .0619 & {\bf.0522} && .0995 & .0784 & .0611 & .0552\\ \cline{1-5} \cline{7-10}
  ATS &     .0640 & .0543 & {\bf.0530} & {\bf.0484} && .0538 & .0462 &.0447& .0458  \\ \cline{1-5} \cline{7-10}
  WTS-Para&  .0798 & .0727 & .0648 & .0604 && .0800 & .0690 & .0638 & .0582 \\ \cline{1-5} \cline{7-10}
  WTS-Wild& .1167 & .0975 & .0786 & .0689 && .1300 & .1083 & .0870 & .0707  \\ \cline{1-5} \cline{7-10}
  WTS-$\chi_{15}^2$   & .5231 & .2387 & .1100 & .0812 && .5448 & .2389 & .1085 & .0764\\ \cline{1-5} \cline{7-10}
  MATS-Para &.0676 & .0622 & .0576 & {\bf.0520} && .0647 & .0579 & .0540 & .0545 \\ \cline{1-5} \cline{7-10}
  MATS-Wild& .0958 & .0754 & .0640 & .0544 && .1036 & .0816 & .0655 & .0605 \\ \cline{1-5} \cline{7-10}
  Bartlett-S &.0166 & .0400 & {\bf.0528} & {\bf.0515} && .0264 & .0594 & .0655 & .0613 \\ \cline{1-5} \cline{7-10}
  Bartlett-P &.0254 & .0361 & .0452 & {\bf.0480} && .0299 & .0405 & .0451 & {\bf.0485}  \\ \cline{1-5} \cline{7-10}
  Box's M-$\chi_{15}^2$&  .1075 & .0976 & .0956 & .0938 && .2707 & .2896 & .3156 & .3250 \\ \cline{1-5} \cline{7-10}
  Box's M-F& .1012 & .0961 & .0952 & .0938 && .2612 & .2881 & .3153 & .3249 \\ \cline{1-5} \cline{7-10}
  \hline
\end{tabular}
\caption{Simulated type-I-error rates ($\alpha=5\%$) in scenario $A)$ ($\mathcal{H}_0^{\vv}:\vV_{1}=\vV_{2}$)  for ATS, WTS, MATS, Bartletts test and Box's M-test.  The observation vectors have dimension 5, covariance matrix $(\vV)_{ij}=0.6^{|i-j|}$ and there is always the same relation between group samples size with $n_1:=0.6\cdot N$ resp. $n_2:=0.4\cdot N$.}\label{tab:SimAAppendix1}
\end{table}
\renewcommand{\baselinestretch}{2}
 For the additional covariance matrix, again the ATS with parametric bootstrap had the best type-I-error control in nearly every setting. Moreover, the influence of the used covariance matrix could be seen, but it neither seemed to be strong nor had a systematical effect on the quality. For the additional distribution $\varphi_{ATS}^*$ exhibited a good performance in most cases, particularly for scenario A). So \Cref{tab:SimAAppendix1}-\Cref{tab:SimBAppendix2} in total confirmed the results from Section 5. The usage and the performance of the MATS showed the variety of our approach one more time.}

\renewcommand{\baselinestretch}{1}
\begin{table}[htbp]
\centering
\begin{tabular}{|l||c|c|c|c||c||c|c|c|c||}\hline
    &\multicolumn{4}{|c||}{$t_9$}&&\multicolumn{4}{|c||}{Normal}\\\cline{1-5} \cline{7-10}
   \hspace*{.1cm}N&   50&100&250&500&&   50&100&250&500
       \\\cline{1-5} \cline{7-10}
ATS-Para & {\bf .0522} & .0544 & {\bf.0517} &{\bf .0504} && .0613 & .0561 & .0533 & .0535  \\ \cline{1-5} \cline{7-10}
  ATS-Wild&  .0772 & .0683 & .0578 & {.0544} && .0782 & .0648 & .0563 & .0541 \\ \cline{1-5} \cline{7-10}
  ATS & .0573 & .0559 & {\bf.0514 }& {\bf.0499} && .0658 & .0575 & .0534 & {\bf.0522} \\ \cline{1-5} \cline{7-10}
  WTS-Para&.0618 & .0641 & .0599 & {.0550} && .0664 & .0665 & .0622 & .0562 \\ \cline{1-5} \cline{7-10}
  WTS-Wild&  .0980 & .0895 & .0726 & .0643 && .0961 & .0852 & .0706 & .0612 \\ \cline{1-5} \cline{7-10}

  WTS-$\chi_{15}^2$ & .4965 & .2168 & .1002 & .0738 && .5000 & .2161 & .1054 & .0757 \\ \cline{1-5} \cline{7-10}
  MATS-Para & .0608 & .0611 &.0553 & .0537 && .0669 & .0603 & .0554 & .0535  \\ \cline{1-5} \cline{7-10}
  MATS-Wild& .0786 & .0699 & .0599 & .0560 && .0837 & .0668 & .0583 & .0553 \\ \cline{1-5} \cline{7-10}
  Bartlett-S &.0171 & {\bf.0488} & .0576 & {\bf.0526} && .0112 & .0368 &{\bf .0481} &{\bf .0482} \\ \cline{1-5} \cline{7-10}
  Bartlett-P &.0233 & .0392 & .0464 & .0465 && .0199 & .0360 & .0452 & .0467 \\ \cline{1-5} \cline{7-10}
  Box's M-$\chi_{15}^2$&.1308 & .1337 & .1361 & .1401& & .0638 & .0575 & {\bf.0521} & {\bf .0496}  \\ \cline{1-5} \cline{7-10}
  Box's M-F&  .1238 & .1322 & .1358 & .1400 && .0609 & .0567 & {\bf.0520} & {\bf.0496} \\ \cline{1-5} \cline{7-10}
   \hline\multicolumn{10}{|c|}{}\\\hline
   &\multicolumn{4}{|c||}{Skew Normal}&&\multicolumn{4}{|c||}{Gamma}\\\cline{1-5} \cline{7-10}
   \hspace*{.1cm}N&   50&100&250&500&&   50&100&250&500
       \\\cline{1-5} \cline{7-10}
ATS-Para &.0602 & {.0543} & .0545 & {\bf.0502} &&{\bf .0502} & {\bf.0475} & {\bf.0473} & {\bf.0484} \\   \cline{1-5} \cline{7-10}
  ATS-Wild&  .0872 & .0687 & .0595 & {\bf.0521} && .0962 & .0749 & .0614 & .0565  \\  \cline{1-5} \cline{7-10}
  ATS &  .0655 & {.0552} &{ .0537} & {\bf.0495} &&{ .0554} &{\bf .0490} & {.0469} &{\bf .0480}\\   \cline{1-5} \cline{7-10}
  WTS-Para&  .0797 & .0729 & .0648 & .0603 && .0813 & .0693 & .0637 & .0580 \\   \cline{1-5} \cline{7-10}
  WTS-Wild& .1167 & .0975 & .0786 & .0689 && .1300 & .1083 & .0870 & .0707 \\   \cline{1-5} \cline{7-10}
  WTS-$\chi_{15}^2$  &.5231 & .2387 & .1100 & .0812 && .5448 & .2389 & .1085 & .0764  \\   \cline{1-5} \cline{7-10}
  MATS-Para & .0689 & .0631 & .0585 & {\bf.0524} && .0675 & .0611 & .0567 & {.0543}  \\   \cline{1-5} \cline{7-10}
  MATS-Wild&  .0889 & .0730 & .0624 &{ .0538} && .0976 & .0787 & .0665 & .0584  \\   \cline{1-5} \cline{7-10}
  Bartlett-S &  .0164 & .0402 & {\bf.0528} & {\bf.0516} && .0264 & .0595 & .0663 & .0612 \\   \cline{1-5} \cline{7-10}
  Bartlett-P &.0254 & .0361 & {.0452} & {\bf.0480} && .0299 & .0405 & {.0451 }& {\bf.0485} \\   \cline{1-5} \cline{7-10}
  Box's M-$\chi_{15}^2$&  .1075 & .0976 & .0956 & .0938 && .2707 & .2896 & .3156 & .3250 \\   \cline{1-5} \cline{7-10}
  Box's M-F&   .1012 & .0961 & .0952 & .0938 && .2612 & .2881 & .3153 & .3249\\   \cline{1-5} \cline{7-10}
  \hline
\end{tabular}
  \caption{Simulated type-I-error rates ($\alpha=5\%$) in scenario $A)$ ($\mathcal{H}_0^{\vv}:\vV_{1}=\vV_{2}$)  for ATS, WTS, MATS, Bartletts test and Box's M-test.  The observation vectors have dimension 5, covariance matrix  $\vV=\vI_5+\vJ_5$ and there is always the same relation between group samples size with $n_1:=0.6\cdot N$ resp. $n_2:=0.4\cdot N$.}
  \label{tab:SimAAppendix2}
\end{table}
\renewcommand{\baselinestretch}{2}
   As expected, it can be seen in \Cref{tab:SimFDiss1} and \Cref{tab:SimFDiss2} that all tests performed generally worse than for just two groups, although some individual results were better. 
In particular, for both Barlett tests and all WTS tests, there was a significant worsening. 
In part, the error rate was almost halved for the Bartlett tests and doubled for the WTS. In comparison, the worsening of Box`s M-test for normal distribution and $\varphi_{ATS}^*$ and $\varphi_{ATS}$, in general, was substantially less pronounced. 
In fact, these were the only tests with error rates in our $95\%$ binomial interval.
So our tests also performed well for this hypothesis although for some distributions bigger sample sizes were required. \\

\renewcommand{\baselinestretch}{1}
    
    \begin{table}[htbp]
\centering
\begin{tabular}{|l||c|c|c|c||c||c|c|c|c||}\hline
    &\multicolumn{4}{|c||}{$t_9$}&&\multicolumn{4}{|c||}{Normal}\\\cline{1-5} \cline{7-10}
   \hspace*{.1cm}N&   25&50&125&250&&   25&50&125&250
       \\\cline{1-5} \cline{7-10}
ATS-Para& .0363 & .0394 & .0395 & .0420 && {.0465} & {\bf.0473} & {\bf.0495} & {\bf.0505} \\ \cline{1-5} \cline{7-10}
  ATS-Wild& .0607 & {.0548} & {\bf.0493} & {\bf.0492} && .0682 & .0573 & {.0542} & {\bf.0527} \\ \cline{1-5} \cline{7-10}
  ATS & .0437 & .0408 & .0403 & .0413 && {.0547} & {\bf.0501} & {\bf.0492} & {\bf.0501} \\ \cline{1-5} \cline{7-10}
  WTS-Para&  .0879 & .0755 & .0636 & .0586 && .0855 & .0702 & .0622 & {.0545}  \\ \cline{1-5} \cline{7-10}
  WTS-Wild& .1193 & .0939 & .0726 & .0653 && .1052 & .0795 & .0660 & {.0557} \\ \cline{1-5} \cline{7-10}
  WTS-$\chi_{4}^2$& .1863 & .1141 & .0751 & .0652 && .1826 & .1109 & .0682 & .0594 \\ \cline{1-5} \cline{7-10}
  MATS-Para &.0832 & .0752 & .0626 & .0576 && .0803 & .0673 & .0585 & {.0544} \\ \cline{1-5} \cline{7-10}
  MATS-Wild& .1251 & .0975 & .0733 & .0643 && .1092 & .0830 & .0639 & .0575  \\\cline{1-5} \cline{7-10}

   \hline\multicolumn{10}{|c|}{}\\\hline
   &\multicolumn{4}{|c||}{Skew Normal}&&\multicolumn{4}{|c||}{Gamma}\\\cline{1-5} \cline{7-10}
   \hspace*{.1cm}N&   25&50&125&250&&   25&50&125&250
       \\\cline{1-5} \cline{7-10}
ATS-Para & {\bf.0481} & .0419 & {.0454} & {\bf.0483} && .0388 & .0363 & .0371 & .0407 \\ \cline{1-5} \cline{7-10}
  ATS-Wild&.0787 & .0618 & {.0550} &{ .0547} && .0805 & .0645 & {.0535} &{\bf .0524}\\ \cline{1-5} \cline{7-10}
  ATS & {.0566} & {.0451} & {.0458} & {\bf.0487} && {.0455} & .0383 & .0373 & .0397 \\ \cline{1-5} \cline{7-10}
  WTS-Para&.1099 & .0886 & .0726 & .0618 && .1441 & .1136 & .0839 & .0711 \\ \cline{1-5} \cline{7-10}
  WTS-Wild& .1458 & .1112 & .0826 & .0675 && .2099 & .1590 & .1076 & .0847 \\ \cline{1-5} \cline{7-10}
  WTS-$\chi_{4}^2$& .2207 & .1277 & .0797 & .0708 && .2609 & .1628 & .0939 & .0761\\ \cline{1-5} \cline{7-10}
  MATS-Para &.0948 & .0756 & .0652 & .0596 && .1046 & .0904 & .0748 & .0674 \\ \cline{1-5} \cline{7-10}
  MATS-Wild& .1431 & .1012 & .0783 & .0655 && .1803 & .1384 & .0960 & .0812\\ \cline{1-5} \cline{7-10}
  \hline
\end{tabular}
 \caption{Simulated type-I-error rates ($\alpha=5\%$) in scenario $B)$ ($\mathcal{H}_0^{\vv}:\vV_{1 11}=\vV_{2 11}=...=\vV_{1 55}$) for ATS, WTS and MATS with 5-dimensional vectors and  $(\vV)_{ij}=0.6^{|i-j|}$.  }  \label{tab:SimBAppendix1}
\end{table}

\begin{table}[htbp]
\centering
\begin{tabular}{|l||c|c|c|c||c||c|c|c|c||}\hline
    &\multicolumn{4}{|c||}{$t_9$}&&\multicolumn{4}{|c||}{Normal}\\\cline{1-5} \cline{7-10}
   \hspace*{.1cm}N&   25&50&125&250&&   25&50&125&250
       \\\cline{1-5} \cline{7-10}
ATS-Para &  .0337 & .0344 & .0383 & {.0429} && {.0446} &{ .0436} & {\bf.0490} & {\bf.0477}  \\   \cline{1-5} \cline{7-10}
  ATS-Wild&.0599 & {.0548} &{\bf .0494} &{\bf .0497} && .0662 &.0563 & .0561 &{\bf .0523}\\   \cline{1-5} \cline{7-10}
  ATS & .0411 & .0381 & .0390 & .0429 && {\bf.0516} & {.0464 }& {\bf.0503} & {\bf.0482}  \\   \cline{1-5} \cline{7-10}
  WTS-Para&  .0885 & .0756 & .0633 & .0602 && .0802 & .0659 & .0603 &{ .0540}  \\   \cline{1-5} \cline{7-10}
  WTS-Wild& .1201 & .0951 & .0740 & .0659 && .1024 & .0760 & .0651 & .0565  \\   \cline{1-5} \cline{7-10}
  WTS-$\chi_{4}^2$ & .1893 & .1123 & .0746 & .0655 && .1794 & .1052 & .0665 & .0589 \\   \cline{1-5} \cline{7-10}
  MATS-Para &  .0816 & .0716 & .0619 & .0585 && .0787 & .0633 & .0588 & {.0533} \\   \cline{1-5} \cline{7-10}
  MATS-Wild& .1330 & .1007 & .0758 & .0662 && .1144 & .0826 & .0666 & .0579\\   \cline{1-5} \cline{7-10}
   \hline\multicolumn{10}{|c|}{}\\\hline
   &\multicolumn{4}{|c||}{Skew Normal}&&\multicolumn{4}{|c||}{Gamma}\\\cline{1-5} \cline{7-10}
   \hspace*{.1cm}N&   25&50&125&250&&   25&50&125&250
       \\\cline{1-5} \cline{7-10}
ATS-Para & {.0446} & .0401 & {.0466} & {\bf.0474} && .0351 & .0325 & .0347 & .0363 \\   \cline{1-5} \cline{7-10}
  ATS-Wild& .0800 & .0617 & .0585 & {.0533} && .0783 & .0634 & {\bf.0523} & {\bf.0488}\\   \cline{1-5} \cline{7-10}
  ATS & {.0544} & .0426 & {.0464} & {\bf.0474} && .0420 & .0343 & .0351 & .0365   \\   \cline{1-5} \cline{7-10}
  WTS-Para&.1110 & .0907 & .0732 & .0627 && .1536 & .1173 & .0847 & .0709\\   \cline{1-5} \cline{7-10}
  WTS-Wild& .1463 & .1136 & .0839 & .0691 && .2234 & .1647 & .1088 & .0843 \\   \cline{1-5} \cline{7-10}
  WTS-$\chi_{4}^2$ &  .2204 & .1301 & .0798 & .0688 && .2720 & .1642 & .0983 & .0771 \\   \cline{1-5} \cline{7-10}
  MATS-Para &.0924 & .0764 & .0680 & .0592 && .1023 & .0877 & .0745 & .0665\\   \cline{1-5} \cline{7-10}
  MATS-Wild& .1487 & .1079 & .0826 & .0660 && .1888 & .1446 & .0994 & .0824  \\   \cline{1-5} \cline{7-10}
  \hline
\end{tabular}
  \caption{Simulated type-I-error rates ($\alpha=5\%$) in  scenario $B)$ ($\mathcal{H}_0^{\vv}:\vV_{1 11}=\vV_{1 22}=...=\vV_{1 55}$) for ATS, WTS and MATS with 5-dimensional vectors and $\vV=\vI_5+\vJ_5$.  }
  \label{tab:SimBAppendix2}
\end{table}

\begin{table}[htbp]
\centering
\begin{tabular}{|l||c|c|c|c||c||c|c|c|c||}\hline
    &\multicolumn{4}{|c||}{$t_9$}&&\multicolumn{4}{|c||}{Normal}\\\cline{1-5} \cline{7-10}
   \hspace*{.1cm}N&   50&100&250&500&&   50&100&250&500
       \\\cline{1-5} \cline{7-10}
ATS-Para & .0645 & {.0556 }&{ .0542} & {\bf.0526} && .0651 & .0581 & {.0537} & {.0539} \\ \cline{1-5} \cline{7-10}
  ATS-Wild& .0678 & .0573 & {.0556} & {.0532 }&& .0686 & .0598 & {.0542} &{ .0545 } \\ \cline{1-5} \cline{7-10}
  ATS & .0739 & .0583 &{ .0546} & {\bf.0521} && .0739 & .0609 & {.0544 }&{ .0541}\\ \cline{1-5} \cline{7-10}
  WTS-Para& .0645 & {.0556} & {.0542 }& {\bf.0526} && .0651 & .0581 & {.0537} & {.0539} \\ \cline{1-5} \cline{7-10}
  WTS-Wild& .0678 & .0573 & {.0556} & {.0532} && .0686 & .0598 & {.0542} &{ .0545 }\\ \cline{1-5} \cline{7-10}
  WTS-$\chi_{1}^2$ & .0734 & .0579 & {.0547 }& {\bf.0519} && .0736 & .0605 & {.0535} & {.0538} \\ \cline{1-5} \cline{7-10}
  MATS-Para & .0696 & .0582 & {.0549} &{ .0532 }&& .0718 & .0602 & {.0544} &{ .0546 }\\ \cline{1-5} \cline{7-10}
  MATS-Wild& .0759 & .0609 & .0564 & {.0535} && .0754 & .0623 &{ .0551} & {.0551} \\ \cline{1-5} \cline{7-10}

   \hline\multicolumn{10}{|c|}{}\\\hline
   &\multicolumn{4}{|c||}{Skew Normal}&&\multicolumn{4}{|c||}{Gamma}\\\cline{1-5} \cline{7-10}
  \hspace*{.1cm}N&   50&100&250&500&&   50&100&250&500
       \\\cline{1-5} \cline{7-10}
ATS-Para &.0690 & .0589 &{\bf .0530} &{\bf .0514} && .0715 & .0628 & {.0540} & {.0552}\\  \cline{1-5} \cline{7-10}
  ATS-Wild& .0738 & .0621 & {.0540} &{\bf .0521 }&& .0848 & .0688 & {.0550 }& {.0556 } \\ \cline{1-5} \cline{7-10}
  ATS & .0779 & .0623 & {.0540} &{\bf .0517} && .0814 & .0655 & {.0540} & {.0538} \\ \cline{1-5} \cline{7-10}
  WTS-Para& .0690 & .0589 & {\bf.0530 }& {\bf.0514} && .0715 & .0628 & {.0540} &{ .0552} \\ \cline{1-5} \cline{7-10}
  WTS-Wild& .0738 & .0621 & {.0540 }&{\bf .0521} && .0848 & .0688 & {.0550} &{ .0556}\\ \cline{1-5} \cline{7-10}
  WTS-$\chi_{1}^2$&  .0775 & .0619 & {.0538} &{\bf .0518} && .0811 & .0651 &{ .0540} & {.0540} \\ \cline{1-5} \cline{7-10}
  MATS-Para &.0752 & .0617 & {.0538} & {\bf.0514 }&& .0782 & .0647 & {.0544 }& {.0552 }\\ \cline{1-5} \cline{7-10}
  MATS-Wild& .0810 & .0649 &{ .0550 }&{\bf .0527} && .0926 & .0718 & .0564 & .0561  \\ \cline{1-5} \cline{7-10}
  \hline
\end{tabular}
\caption{Simulated type-I-error rates ($\alpha=5\%$) in scenario $C)$ ($\mathcal{H}_0^{\vv}:\tr(\vV_{1})=\tr(\vV_{2})$)  for ATS, WTS, and MATS.  The observation vectors have dimension 5, covariance matrix $(\vV)_{ij}=0.6^{|i-j|}$ and there is always the same relation between group samples size with $n_1:=0.6\cdot N$ resp. $n_2:=0.4\cdot N$.}\label{tab:SimCAppendix1}
\end{table}

\begin{table}[htbp]
\centering
\begin{tabular}{|l||c|c|c|c||c||c|c|c|c||}\hline
    &\multicolumn{4}{|c||}{$t_9$}&&\multicolumn{4}{|c||}{Normal}\\\cline{1-5} \cline{7-10}
   \hspace*{.1cm}N&   50&100&250&500&&   50&100&250&500
       \\\cline{1-5} \cline{7-10}
ATS-Para & .0607 & {.0539} & {.0545} & {.0532} && .0632 & .0584 & {.0545} &{ .0555}  \\\cline{1-5} \cline{7-10}
  ATS-Wild& .0666 & .0570 & {.0560} &{ .0542} && .0674 & .0609 & {.0544} & {.0548} \\ \cline{1-5} \cline{7-10}
  ATS & .0697 & .0573 & {.0555} & {\bf.0527 }&& .0724 & .0620 & {.0539 }&{ .0550 } \\ \cline{1-5} \cline{7-10}
  WTS-Para& .0607 & .0539 & {.0545} &{ .0532} && .0632 & .0584 & {.0545 }& {.0555}  \\ \cline{1-5} \cline{7-10}
  WTS-Wild& .0666 & .0570 & {.0560} & {.0542 }&& .0674 & .0609 & {.0544} & {.0548 } \\ \cline{1-5} \cline{7-10}
  WTS-$\chi_{1}^2$ &.0698 & .0575 & {.0551 }& {\bf.0527} && .0720 & .0619 & {.0542} &{ .0546 }\\ \cline{1-5} \cline{7-10}
  MATS-Para& .0668 & .0569 & .0563 & {.0537} && .0687 & .0613 &{ .0554} & .0561 \\ \cline{1-5} \cline{7-10}
  MATS-Wild& .0742 & .0599 & .0569 & {.0550} && .0741 & .0635 & {.0555} &{ .0554} \\ \cline{1-5} \cline{7-10}

   \hline\multicolumn{10}{|c|}{}\\\hline
   &\multicolumn{4}{|c||}{Skew Normal}&&\multicolumn{4}{|c||}{Gamma}\\\cline{1-5} \cline{7-10}
   \hspace*{.1cm}N&   50&100&250&500&&   50&100&250&500
       \\\cline{1-5} \cline{7-10}
ATS-Para & .0685 & .0575 & {.0533} &{\bf .0521} && .0713 & .0603 & {\bf.0528} & {.0540} \\ \cline{1-5} \cline{7-10}
  ATS-Wild&  .0747 & .0620 & {.0536} & {\bf.0530} && .0828 & .0677 & .0562 & {.0558} \\ \cline{1-5} \cline{7-10}
  ATS &  .0765 & .0607 & {.0536} & {\bf.0518} && .0804 & .0631 &{ .0541} &{ .0544}\\ \cline{1-5} \cline{7-10}
  WTS-Para&.0685 & .0575 &{ .0533} & {\bf.0521} && .0713 & .0603 &{\bf .0528} & {.0540 }\\ \cline{1-5} \cline{7-10}
  WTS-Wild&  .0747 & .0620 & {.0536 }& {\bf.0530} && .0828 & .0677 & .0562 & {.0558} \\ \cline{1-5} \cline{7-10}
  WTS-$\chi_{1}^2$  &.0766 & .0612 &{ .0536} &{\bf .0518} && .0804 & .0630 & {.0537} &{ .0543}  \\ \cline{1-5} \cline{7-10}
  MATS-Para & .0742 & .0608 &{ .0543} & {\bf.0524 }&& .0773 & .0627 & {.0540} &{ .0547}  \\ \cline{1-5} \cline{7-10}
  MATS-Wild& .0812 & .0645 & {.0547} &{ .0536} && .0923 & .0708 & .0569 & .0563  \\ \cline{1-5} \cline{7-10}
  \hline
\end{tabular}
\caption{Simulated type-I-error rates ($\alpha=5\%$) in scenario $C)$ ($\mathcal{H}_0^{\vv}:\tr(\vV_{1})=\tr(\vV_{2})$)  for ATS, WTS, and MATS.  The observation vectors have dimension 5, covariance matrix $\vV=\vI_5+\vJ_5$ and there is always the same relation between group samples size with $n_1:=0.6\cdot N$ resp. $n_2:=0.4\cdot N$.}\label{tab:SimCAppendix2}
\end{table}

\begin{table}[tbp]
\centering
\begin{tabular}{|l||c|c|c|c||c||c|c|c|c||}\hline
    &\multicolumn{4}{|c||}{$t_9$}&&\multicolumn{4}{|c||}{Normal}\\\cline{1-5} \cline{7-10}
   \hspace*{.1cm}N&   80&160&400&800&&  80&160&400&800
       \\\cline{1-5} \cline{7-10}
ATS-Para & .0415 & .0409 & .0435 & .0467 && .0531 & {\bf .0511} &{\bf  .0487} & {\bf .0505} \\ \cline{1-5} \cline{7-10}
  ATS-Wild & .0776 & .0643 & .0562 & .0544 && .0809 & .0666 & .0547 & .0536 \\ \cline{1-5} \cline{7-10}
  ATS & .0452 & .0422 & .0436 & .0461 && .0578 & {\bf .0527} & {\bf .0474 }& {\bf .0501 }\\ \cline{1-5} \cline{7-10}
  WTS-Para & .1186 & .1083 & .0831 & .0703 && .1341 & .1198 & .0853 & .0685 \\ \cline{1-5} \cline{7-10}
  WTS-Wild & .1821 & .1533 & .1031 & .0818 && .1869 & .1514 & .0976 & .0743 \\ \cline{1-5} \cline{7-10}
  WTS-$\chi_{30}^2$ & .8619 & .4602 & .1964 & .1258 && .8648 & .4556 & .1893 & .1239 \\ \cline{1-5} \cline{7-10}
  MATS-Para & .0763 & .0673 & .0580 & .0553 && .0825 & .0689 & .0573 & .0544 \\ \cline{1-5} \cline{7-10}
  MATS-Wild & .1212 & .0871 & .0687 & .0596 && .1165 & .0855 & .0641 & .0574 \\ \cline{1-5} \cline{7-10}
  Bartlett-S & .0117 & .0199 & .0299 & .0379 && .0204 & .0301 & .0401 & .0432 \\ \cline{1-5} \cline{7-10}
  Bartlett-P & .0226 & .0309 & .0371 & .0444 && .0321 & .0396 & .0466 & .0465 \\ \cline{1-5} \cline{7-10}
  Box's M-$\chi_{30}^2$ & .1580 & .1678 & .1825 & .1872 && .0671 & .0560 &{\bf  .0497} &{\bf  .0517} \\ \cline{1-5} \cline{7-10}
  Box's M-F & .1521 & .1653 & .1824 & .1872 && .0638 & .0546 & {\bf .0496} &{\bf  .0517} \\ 

   \hline\multicolumn{10}{|c|}{}\\\hline
   &\multicolumn{4}{|c||}{Skew Normal}&&\multicolumn{4}{|c||}{Gamma}\\\cline{1-5} \cline{7-10}
  \hspace*{.1cm}N&   80&160&400&800&&   80&160&400&800
       \\\cline{1-5} \cline{7-10}
ATS-Para & .0454 & .0459 & .0465 &{\bf  .0482 }&& .0349 & .0328 & .0373 & .0425 \\ \cline{1-5} \cline{7-10}
  ATS-Wild & .0858 & .0704 & .0575 & .0548 && .0938 & .0716 & .0574 & .0564 \\ \cline{1-5} \cline{7-10}
  ATS &{\bf  .0505} & {\bf .0474} & .0465 & {\bf .0480} && .0401 & .0335 & .0372 & .0425 \\ \cline{1-5} \cline{7-10}
  WTS-Para & .1461 & .1332 & .0933 & .0748 && .1531 & .1420 & .1039 & .0830 \\ \cline{1-5} \cline{7-10}
  WTS-Wild & .2047 & .1745 & .1116 & .0866 && .2331 & .2060 & .1391 & .1024 \\ \cline{1-5} \cline{7-10}
  WTS-$\chi_{30}^2$ & .8823 & .4911 & .2031 & .1309 && .8992 & .5254 & .2220 & .1417 \\ \cline{1-5} \cline{7-10}
  MATS-Para & .0873 & .0727 & .0586 & .0582 && .0913 & .0795 & .0658 & .0605 \\ \cline{1-5} \cline{7-10}
  MATS-Wild & .1322 & .0972 & .0686 & .0644 && .1559 & .1179 & .0836 & .0713 \\ \cline{1-5} \cline{7-10}
  Bartlett-S & .0139 & .0224 & .0343 & .0402 && .0066 & .0110 & .0239 & .0319 \\ \cline{1-5} \cline{7-10}
  Bartlett-P & .0268 & .0326 & .0409 & .0448 && .0161 & .0224 & .0340 & .0377 \\ \cline{1-5} \cline{7-10}
  Box's M-$\chi_{30}^2$ & .1211 & .1156 & .1138 & .1160 && .3819 & .4306 & .4661 & .4776 \\ \cline{1-5} \cline{7-10}
  Box's M-F & .1148 & .1139 & .1136 & .1159 && .3723 & .4281 & .4657 & .4775 \\ 
  \hline
\end{tabular}
\caption{Simulated type-I-error rates ($\alpha=5\%$) in scenario $E)$ ($\mathcal{H}_0^{\vv}:\vV_{1}=\vV_{2}=\vV_3$)  for ATS, WTS, and MATS.  The observation vectors have dimension 5, covariance matrix $(\vV)_{ij}=0.6^{|i-j|}$ and there is always the same relation between group samples size with $n_1:=0.4\cdot N$, $n_2:=0.25\cdot N$ and $n_3:=0.35\cdot N$.}\label{tab:SimFDiss1}
\end{table}

\begin{table}[tbp]
\centering
\begin{tabular}{|l||c|c|c|c||c||c|c|c|c||}\hline
    &\multicolumn{4}{|c||}{$t_9$}&&\multicolumn{4}{|c||}{Normal}\\\cline{1-5} \cline{7-10}
   \hspace*{.1cm}N&   80&160&400&800&&  80&160&400&800
       \\\cline{1-5} \cline{7-10}
ATS-Para & .0435 & .0435 & .0456 &{\bf  .0474 }&& .0561 & {\bf .0521} & {\bf .0480 }& {\bf .0518} \\ \cline{1-5} \cline{7-10}
  ATS-Wild & .0749 & .0618 & .0552 & {\bf .0530} && .0782 & .0639 & .0533 & .0537 \\ \cline{1-5} \cline{7-10}
  ATS & .0464 & .0444 & .0452 & .0468 && .0606 & .0536 & {\bf .0482} & {\bf .0505} \\ \cline{1-5} \cline{7-10}
  WTS-Para & .1190 & .1089 & .0842 & .0693 && .1340 & .1190 & .0845 & .0693 \\ \cline{1-5} \cline{7-10}
  WTS-Wild & .1821 & .1533 & .1031 & .0818 && .1869 & .1514 & .0976 & .0743 \\ \cline{1-5} \cline{7-10}
  WTS-$\chi_{30}^2$ & .8619 & .4602 & .1964 & .1258 && .8648 & .4556 & .1893 & .1239 \\ \cline{1-5} \cline{7-10}
  MATS-Para & .0821 & .0695 & .0575 & .0555 && .0853 & .0726 & .0583 & .0557 \\ \cline{1-5} \cline{7-10}
  MATS-Wild & .1169 & .0856 & .0650 & .0585 && .1133 & .0857 & .0628 & .0585 \\ \cline{1-5} \cline{7-10}
  Bartlett-S & .0117 & .0199 & .0299 & .0379 && .0204 & .0301 & .0401 & .0432 \\ \cline{1-5} \cline{7-10}
  Bartlett-P & .0226 & .0309 & .0371 & .0444 && .0321 & .0396 & .0466 & .0465 \\ \cline{1-5} \cline{7-10}
  Box's M-$\chi_{30}^2$ & .1580 & .1678 & .1825 & .1872 && .0671 & .0560 &{\bf  .0497} & {\bf .0517} \\ \cline{1-5} \cline{7-10}
  Box's M-F & .1521 & .1653 & .1824 & .1872 && .0638 & .0546 & {\bf .0496} & {\bf .0517} \\\cline{1-5} \cline{7-10}

   \hline\multicolumn{10}{|c|}{}\\\hline
   &\multicolumn{4}{|c||}{Skew Normal}&&\multicolumn{4}{|c||}{Gamma}\\\cline{1-5} \cline{7-10}
  \hspace*{.1cm}N&   80&160&400&800&&  80&160&400&800
       \\\cline{1-5} \cline{7-10}
ATS-Para & {\bf .0487} &{\bf  .0467 }&{\bf  .0475 }&{\bf  .0494} && .0379 & .0361 & .0413 & .0455 \\ \cline{1-5} \cline{7-10}
  ATS-Wild & .0809 & .0660 & .0566 & .0560 && .0901 & .0698 & .0599 & .0559 \\ \cline{1-5} \cline{7-10}
  ATS & .0535 & {\bf .0487} & {\bf .0481 }& {\bf .0503} && .0401 & .0363 & .0404 & .0456 \\ \cline{1-5} \cline{7-10}
  WTS-Para & .1461 & .1325 & .0926 & .0749 && .1539 & .1431 & .1054 & .0827 \\ \cline{1-5} \cline{7-10}
  WTS-Wild & .2047 & .1745 & .1116 & .0866 && .2331 & .2060 & .1391 & .1024 \\ \cline{1-5} \cline{7-10}
  WTS-$\chi_{30}^2$ & .8823 & .4911 & .2031 & .1309 && .8992 & .5254 & .2220 & .1417 \\ \cline{1-5} \cline{7-10}
  MATS-Para & .0891 & .0740 & .0600 & .0580 && .0947 & .0818 & .0676 & .0603 \\ \cline{1-5} \cline{7-10}
  MATS-Wild & .1244 & .0907 & .0686 & .0634 && .1456 & .1123 & .0810 & .0690 \\ \cline{1-5} \cline{7-10}
  Bartlett-S & .0139 & .0224 & .0343 & .0402 && .0066 & .0110 & .0239 & .0319 \\ \cline{1-5} \cline{7-10}
  Bartlett-P & .0268 & .0326 & .0409 & .0448 && .0161 & .0224 & .0340 & .0377 \\ \cline{1-5} \cline{7-10}
  Box's M-$\chi_{30}^2$ & .1211 & .1156 & .1138 & .1160 && .3819 & .4306 & .4661 & .4776 \\ \cline{1-5} \cline{7-10}
  Box's M-F & .1148 & .1139 & .1136 & .1159 && .3723 & .4281 & .4657 & .4775 \\ 
  \hline
\end{tabular}
\caption{Simulated type-I-error rates ($\alpha=5\%$) in scenario $E)$ ($\mathcal{H}_0^{\vv}:\vV_{1}=\vV_{2}=\vV_3$)  for ATS, WTS, and MATS.  The observation vectors have dimension 5, covariance matrix $(\vV)=\vI_5+\vJ_5$ and there is always the same relation between group samples size with $n_1:=0.4\cdot N$, $n_2:=0.25\cdot N$ and $n_3:=0.35\cdot N$.}\label{tab:SimFDiss2}
\end{table}

At last, we investigate whether the performance stayed essentially the same if the dimension was increased as long as the relation between sample size and dimension remained the same.
Therefore we considered dimension $d=7$ which led to $p=28$. Note that this is substantially larger than for $d=5$ where we had $p=15$.
With $\vN=(70,140,350,700)$ we considered the setting from $A)$ ($\mathcal{H}_0^{\vv}:\vV_{1}=\vV_{2}$) with  the same kind of distributions and covariance matrices, but for dimension 7.
The corresponding results are displayed in \Cref{tab:SimA7Diss1}  and \Cref{tab:SimA7Diss2}.

It is interesting, that the higher dimension together with the larger number of observations improved the results of our test for some distributions like for the skew normal distribution, while for others the performance deteriorated, e.g., for the gamma distribution.
But across all distributions, the quality of the tests' performance was comparable to the situation with dimension 5.
In some way, this is surprising because the sample size was increased in linear relation to the dimension $d$, and not in relation to the dimension of the vectorized covariance matrix $p$, where the latter grows much more rapidly.  
Considering the dimension of this vector as the decisive factor, the relative sample sizes are clearly lower than for dimension 5. Once more, this simulation demonstrated good small sample behavior.
By contrast, the results of both Bartlett statistics showed the impact of this smaller sample size in relation to $p$. In particular, for $N=70$ and $N=140$, a worse type-I-error control could be seen, for example, for the skewed normal distribution with autoregressive covariance matrix. Moreover, Box's M test with both kinds of critical values performed considerably worse in the case of normally distributed data. While for larger sample sizes in the case of dimension 5 all error rates were in the $95\%$ binomial interval, for dimension 7 none of them was in the interval.\\

The results from this setting show that in cases of higher dimension, the performance for smaller sample sizes is of essential importance. Unfortunately, most of the existing procedures do not perform in a satisfactory way in this situation. \\\\
To sum up, $\varphi_{ATS}^*$ and $ \varphi_{ATS}$ led to good finite sample results, even for small sample sizes and challenging null hypotheses or higher dimension. 
The excellent small sample approximation and the variety of applicable situations for this approach make the results for the ATS with parametric bootstrap even a little bit more convincing.  The tests from \cite{zhang1993} were inadequate in most of the cases. Also, neither the WTS nor the MATS based tests were reliable choices for small to moderate sample size settings.

In comparison, the ATS with parametric bootstrap as well as based on a Monte-Carlo simulation exhibited rather good results in particular for higher dimension. Moreover, for more groups, these tests were the only ones with a sufficiently convincing performance in case of non-normality. Thus, these additional simulations emphasized again the wide applicability of both of these newly developed tests.

\begin{table}[tbp]
\centering
\begin{tabular}{|l||c|c|c|c||c||c|c|c|c||}\hline
    &\multicolumn{4}{|c||}{$t_9$}&&\multicolumn{4}{|c||}{Normal}\\\cline{1-5} \cline{7-10}
   \hspace*{.1cm}N&   70&140&350&700&&   70&140&350&700
       \\\cline{1-5} \cline{7-10}
ATS-Para & {\bf .0474} & .0467 &{\bf  .0492} & {\bf .0500} && .0562 & .0532 & {\bf .0527} & {\bf .0482} \\ \cline{1-5} \cline{7-10}
  ATS-Wild & .0790 & .0643 & .0576 & .0553 && .0819 & .0637 & .0572 & {\bf .0491} \\ \cline{1-5} \cline{7-10}
  ATS & {\bf .0495} & .0469 & {\bf .0495} & {\bf .0494} && .0589 & .0543 & {\bf .0528} & {\bf .0471} \\ \cline{1-5} \cline{7-10}
  WTS-Para & .0639 & .0623 & .0599 & .0591 && .0664 & .0715 & .0638 & .0589 \\ \cline{1-5} \cline{7-10}
  WTS-Wild & .0922 & .0863 & .0720 & .0668 && .0937 & .0892 & .0722 & .0642 \\ \cline{1-5} \cline{7-10}
  WTS-$\chi_{28}^2$  & .7933 & .3668 & .1403 & .0906 && .7963 & .3703 & .1428 & .0887 \\ \cline{1-5} \cline{7-10}
  MATS-Para & .0562 & .0534 & .0531 & {\bf .0514} && .0602 & .0565 & {\bf .0530} & {\bf .0484 }\\ \cline{1-5} \cline{7-10}
  MATS-Wild & .0809 & .0680 & .0589 & .0555 && .0828 & .0678 & .0566 & {\bf .0504} \\   \cline{1-5} \cline{7-10}
   Bartlett-S & .0128 & .0392 & {\bf .0508} & .0534 && .0079 & .0323 &{\bf  .0474} & {\bf .0501} \\  \cline{1-5} \cline{7-10}
  Bartlett-P & .0192 & .0323 & .0416 & {\bf .0472} && .0175 & .0355 & .0441 & {\bf .0478} \\  \cline{1-5} \cline{7-10}
  Box's M-$\chi_{28}^2$ & .1440 & .1396 & .1448 & .1442 && .0719 & .0618 & .0562 & .0539 \\  \cline{1-5} \cline{7-10}
  Box's M-F & .1353 & .1374 & .1445 & .1441 && .0654 & .0608 & .0561 & .0539 \\ 

   \hline\multicolumn{10}{|c|}{}\\\hline
   &\multicolumn{4}{|c||}{Skew Normal}&&\multicolumn{4}{|c||}{Gamma}\\\cline{1-5} \cline{7-10}
  \hspace*{.1cm}N&   70&140&350&700&&   70&140&350&700
       \\\cline{1-5} \cline{7-10}
ATS-Para & {\bf .0486} & {\bf .0484} & {\bf .0489} & {\bf .0499} && .0389 & .0406 & .0421 & .0454 \\ \cline{1-5} \cline{7-10}
  ATS-Wild & .0804 & .0654 & .0554 & .0533 && .0881 & .0705 & .0567 & .0542 \\ \cline{1-5} \cline{7-10}
  ATS & {\bf .0508} & {\bf .0494} & {\bf .0476} & {\bf .0487} && .0410 & .0406 & .0417 & .0453 \\ \cline{1-5} \cline{7-10}
  WTS-Para & .0726 & .0758 & .0635 & .0590 && .0747 & .0731 & .0675 & .0606 \\ \cline{1-5} \cline{7-10}
  WTS-Wild & .1054 & .0972 & .0741 & .0659 && .1148 & .1075 & .0854 & .0728 \\ \cline{1-5} \cline{7-10}
  WTS-$\chi_{28}^2$  & .8072 & .3830 & .1431 & .0907 && .8276 & .3935 & .1502 & .0941 \\ \cline{1-5} \cline{7-10}
  MATS-Para & .0583 & .0562 & {\bf .0516} & {\bf .0510} && .0567 & {\bf .0524} & .0541 & {\bf .0497} \\ \cline{1-5} \cline{7-10}
  MATS-Wild & .0858 & .0695 & .0586 & .0546 && .0913 & .0738 & .0637 & .0568 \\ \cline{1-5} \cline{7-10}
    Bartlett-S & .0107 & .0354 & .0467 &{\bf  .0481 }&& .0216 & {\bf .0519} & .0613 & .0590 \\ \cline{1-5} \cline{7-10}
  Bartlett-P & .0184 & .0339 & .0425 & .0459 && .0245 & .0360 & .0448 & .0462 \\ \cline{1-5} \cline{7-10}
  Box's M-$\chi_{28}^2$ & .1123 & .1027 & .0926 & .0939 && .3046 & .3241 & .3417 & .3557 \\ \cline{1-5} \cline{7-10}
  Box's M-F & .1045 & .1012 & .0924 & .0939 && .2933 & .3204 & .3413 & .3555 \\ 
  \hline
\end{tabular}
\caption{Simulated type-I-error rates ($\alpha=5\%$) in scenario $A)$ ($\mathcal{H}_0^{\vv}:\vV_{1}=\vV_{2}$)  for ATS, WTS, and MATS.  The observation vectors have dimension 7, covariance matrix $(\vV)_{ij}=0.6^{|i-j|}$ and there is always the same relation between group samples size with $n_1:=0.6\cdot N$, resp. $n_2:=0.4\cdot N$.}\label{tab:SimA7Diss1}
\end{table}
\begin{table}[tbp]
\centering
\begin{tabular}{|l||c|c|c|c||c||c|c|c|c||}\hline
    &\multicolumn{4}{|c||}{$t_9$}&&\multicolumn{4}{|c||}{Normal}\\\cline{1-5} \cline{7-10}
   \hspace*{.1cm}N&  70&140&350&700&&   70&140&350&700
       \\\cline{1-5} \cline{7-10}
ATS-Para & .0540 & {\bf .0504} & {\bf .0506} & {\bf .0514} && .0622 & .0553 & .0533 & {\bf .0486} \\ \cline{1-5} \cline{7-10}
  ATS-Wild & .0734 & .0611 & .0553 & .0543 && .0757 & .0614 & .0558 & {\bf .0495} \\ \cline{1-5} \cline{7-10}
  ATS & .0564 & {\bf .0511} & {\bf .0505} & {\bf .0512} && .0645 & .0562 & {\bf .0524} & {\bf .0482} \\ \cline{1-5} \cline{7-10}
  WTS-Para & .0626 & .0618 & .0606 & .0597 && .0682 & .0721 & .0642 & .0585 \\ \cline{1-5} \cline{7-10}
  WTS-Wild & .0922 & .0863 & .0720 & .0668 && .0937 & .0892 & .0722 & .0642 \\ \cline{1-5} \cline{7-10}
  WTS-$\chi_{28}^2$  & .7933 & .3668 & .1403 & .0906 && .7963 & .3703 & .1428 & .0887 \\ \cline{1-5} \cline{7-10}
  MATS-Para & .0624 & .0572 & .0535 & {\bf .0526} && .0663 & .0593 & .0549 & {\bf .0493} \\ \cline{1-5} \cline{7-10}
 MATS-Wild & .0760 & .0635 & .0569 & .0548 && .0775 & .0642 & .0564 & {\bf .0497} \\ \cline{1-5} \cline{7-10}
 
   Bartlett-S & .0128 & .0392 & {\bf .0508} & .0534 && .0079 & .0323 & {\bf .0474} &{\bf  .0501} \\  \cline{1-5} \cline{7-10}
  Bartlett-P & .0192 & .0323 & .0416 & {\bf .0472} && .0175 & .0355 & .0441 & {\bf .0478} \\  \cline{1-5} \cline{7-10}
  Box's M-$\chi_{28}^2$ & .1440 & .1396 & .1448 & .1442 && .0719 & .0618 & .0562 & .0539 \\  \cline{1-5} \cline{7-10}
  Box's M-F & .1353 & .1374 & .1445 & .1441 && .0654 & .0608 & .0561 & .0539 \\  

   \hline\multicolumn{10}{|c|}{}\\\hline
   &\multicolumn{4}{|c||}{Skew Normal}&&\multicolumn{4}{|c||}{Gamma}\\\cline{1-5} \cline{7-10}
  \hspace*{.1cm}N&  70&140&350&700&&   70&140&350&700
       \\\cline{1-5} \cline{7-10}
ATS-Para & .0560 & .0546 & {\bf .0507} & {\bf .0494} && .0466 & .0467 & .0467 & {\bf .0483} \\ \cline{1-5} \cline{7-10}
  ATS-Wild & .0746 & .0632 & .0544 & {\bf .0509 }&& .0806 & .0657 & .0556 & {\bf .0527} \\ \cline{1-5} \cline{7-10}
  ATS & .0588 & .0554 & {\bf .0504} & {\bf .0487} && {\bf .0495 }& {\bf .0477} & .0463 & .0467 \\ \cline{1-5} \cline{7-10}
  WTS-Para & .0728 & .0744 & .0628 & .0589 && .0745 & .0726 & .0673 & .0614 \\ \cline{1-5} \cline{7-10}
  WTS-Wild & .1054 & .0972 & .0741 & .0659 && .1148 & .1075 & .0854 & .0728 \\ \cline{1-5} \cline{7-10}
  WTS-$\chi_{28}^2$  & .8072 & .3830 & .1431 & .0907 && .8276 & .3935 & .1502 & .0941 \\ \cline{1-5} \cline{7-10}
  MATS-Para & .0638 & .0595 & .0541 & {\bf .0504} && .0624 & .0592 & .0557 & {\bf .0512} \\ \cline{1-5} \cline{7-10}
  MATS-Wild & .0770 & .0658 & .0566 &{\bf  .0521} && .0828 & .0699 & .0604 & .0550 \\ \cline{1-5} \cline{7-10}
    Bartlett-S & .0107 & .0354 & .0467 & {\bf .0481} && .0216 & {\bf .0519} & .0613 & .0590 \\ \cline{1-5} \cline{7-10}
  Bartlett-P & .0184 & .0339 & .0425 & .0459 && .0245 & .0360 & .0448 & .0462 \\ \cline{1-5} \cline{7-10}
  Box's M-$\chi_{28}^2$ & .1123 & .1027 & .0926 & .0939 && .3046 & .3241 & .3417 & .3557 \\ \cline{1-5} \cline{7-10}
  Box's M-F & .1045 & .1012 & .0924 & .0939 &&.2933 & .3204 & .3413 & .3555 \\ 
  \hline
\end{tabular}
\caption{Simulated type-I-error rates ($\alpha=5\%$) in scenario $A)$ ($\mathcal{H}_0^{\vv}:\vV_{1}=\vV_{2}$)  for ATS, WTS, and MATS.  The observation vectors have dimension 7, covariance matrix $(\vV)=\vI_7+\vJ_7$ and there is always the same relation between group samples size with $n_1:=0.6\cdot N$, resp. $n_2:=0.4\cdot N$.}\label{tab:SimA7Diss2}
\end{table}
\renewcommand{\baselinestretch}{2}

\subsection{More Power Plots}
For the power simulation, we again considered, on the one hand, $\varphi_{ATS}^*$, $\varphi_{ATS}^\star$ and $\varphi_{ATS}$, and on the other hand, $\varphi_{B-S}$ and $\varphi_{B-P}$, all with small ($N=25$ resp. $N=50$) as well as moderate(N=50 resp. N=100) sample size. The results from Section 5 showed that the one-point-alternative is of greater interest, so we just considered this alternative. Besides the skewed normal distribution, we considered the gamma distribution and all hypotheses from Section 5.

For the 100 observations, we see in \Cref{Power1}, that all tests detected the deviation from the null hypothesis earlier, especially for the Bartlett test statistics, which showed bad results for $N=50$. Moreover, here the distance between the different kinds of bootstrap gets smaller for the parametric and the wild bootstrap as well as for the separate and pooled.
In \Cref{Power2} the results for the gamma distribution were similar, but slightly more difference between the kind of bootstraps overall sample sizes.

In \Cref{Power3} and \Cref{Power4} again the wild bootstrap had the best power because of his liberal behavior. While $\varphi_{ATS}$ behaved at the beginning as $\varphi_{ATS}^*$, for larger values of $\delta$ it had obvious more power because it got closer to $\varphi_{ATS}^\star$. The sample size seems to have less impact on this, and the overall the difference in the results for skewed normal distribution and gamma distribution is nearly solely, that the difference between $\varphi_{ATS}^\star$ and $\varphi_{ATS}$ was notable smaller for the skewed normal distribution.

Finally for scenario C) it is really noteworthy that the less liberal test $\varphi_{ATS}^*$ for $\delta>0.5$ had clearly more power than $\varphi_{ATS}^\star$ and $\varphi_{ATS}$. This holds for both sample sizes and both distributions, while it was slightly smaller for the gamma distribution.
Moreover, as for the type-I-error rate, it could be seen that for scenario $C)$ the Monte-Carlo test $\varphi_{ATS}$ was similar to $\varphi_{ATS}^\star$ and not like for the hypothesis of equal covariances similar to $\varphi_{ATS}^*$.\\\\

 It turned out that the power of the ATS is always higher than from Bartlett's statistic unattached from the chosen bootstrap technique. Although in the other hypotheses where Bartlett's test statistic can not be used, the ATS in particular with the wild bootstrap had quite good power curves. The fact that the less liberal $\varphi_{ATS}^*$ had higher power is really interesting.
 Summarizing the results from this extended simulation, our test showed good power even for this hypothesis, which is difficult to detect, and sample sizes, which are really small for dimension $p=15$.

It is important to mention the fact that multiplication with the diagonal matrix changes, not even $\vV$ but also $\vSigma$. So for each $\vDelta$ the eigenvalues of $\vSigma^{1/2}\vC^\top\vC \vSigma^{1/2}$ changed, which consequently changed the limit distribution of the ATS.
Regrettably changing $\vV$ without changing $\vSigma$ is nearly impossible, and there exists no good approach to check the power in situations like this so far. Because of this, for example, in \cite{zhang1993}, just one matrix is used for the calculation of power instead of a whole sequence. Therefore our approach is pretty advanced and kind of intuitive.

\renewcommand{\baselinestretch}{1}
\begin{figure}[htbp]
\begin{minipage}{\textwidth}
     \centering
\includegraphics[trim= 7mm 13mm 5mm 20mm,clip,scale=0.79]{P3OPAAk}
 \end{minipage}
  \begin{minipage}{\textwidth}
     \centering
     \includegraphics[trim= 7mm 12mm 5mm 20mm,clip,scale=0.79]{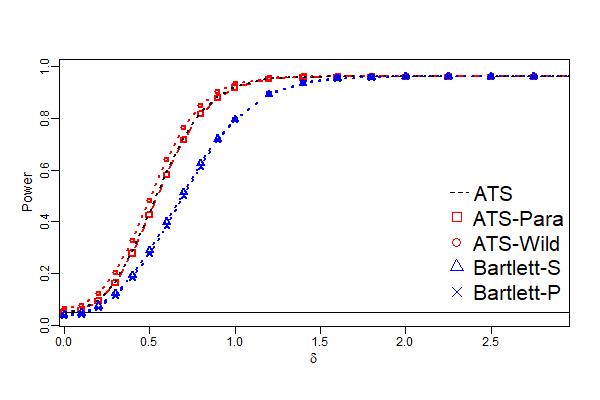}
\end{minipage}

\caption{{Simulated power for an one-point-alternative in scenario $A)$ ($\mathcal{H}_0^{\vv}:\vV_{1}=\vV_{2}$) { for the ATS tests based upon wild bootstrap, parametric bootstrap and Monte-Carlo critical values} as well as the two bootstrap tests based on Bartlett's statistic. The d=5 dimensional error terms are based on the skewed normal distribution with  covariance matrix  $(\vV)_{ij}=0.6^{|i-j|}$ and sample sizes $n_1=30 , n_2=20$ in the first row and $n_1=60, n_2=40$ in the second.}}
\label{Power1}

\end{figure}

\begin{figure}[htbp]
\begin{minipage}{\textwidth}
     \centering
\includegraphics[trim= 7mm 13mm 5mm 20mm,clip,scale=0.79]{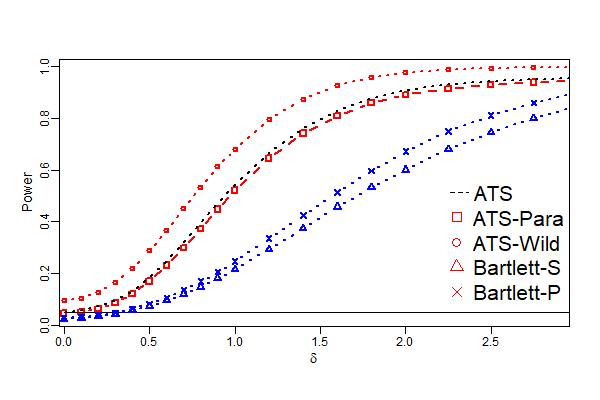}
 \end{minipage}
  \begin{minipage}{\textwidth}
     \centering
      \includegraphics[trim= 7mm 12mm 5mm 20mm,clip,scale=0.79]{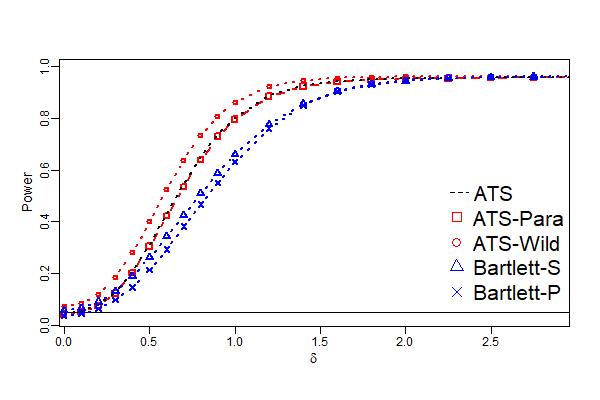}
\end{minipage}

\caption{ Simulated power for a one-point-alternative in scenario $A)$ ($\mathcal{H}_0^{\vv}:\vV_{1}=\vV_{2}$) for the {ATS tests based upon wild bootstrap, parametric bootstrap and Monte-Carlo critical values} as well as the two bootstrap tests based on Bartlett's statistic. The d=5 dimensional error terms are based on the gamma distribution with  covariance matrix  $(\vV)_{ij}=0.6^{|i-j|}$ and sample sizes $n_1=30 , n_2=20$ in the first row and $n_1=60, n_2=40$ in the second.}
\label{Power2}

\end{figure}

\begin{figure}[htbp]
\begin{minipage}{\textwidth}
     \centering
\includegraphics[trim= 7mm 13mm 5mm 20mm,clip,scale=0.79]{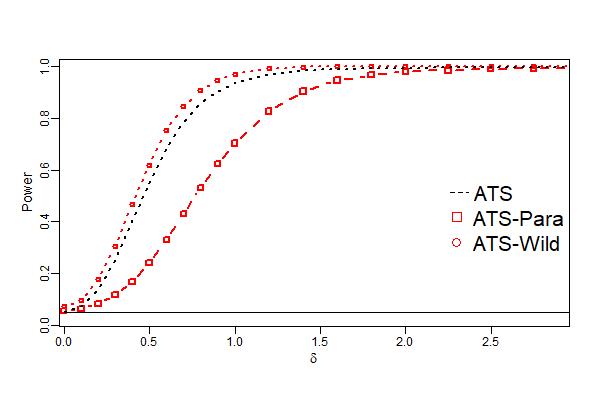}
 \end{minipage}
  \begin{minipage}{\textwidth}
     \centering
           \includegraphics[trim= 7mm 12mm 5mm 20mm,clip,scale=0.79]{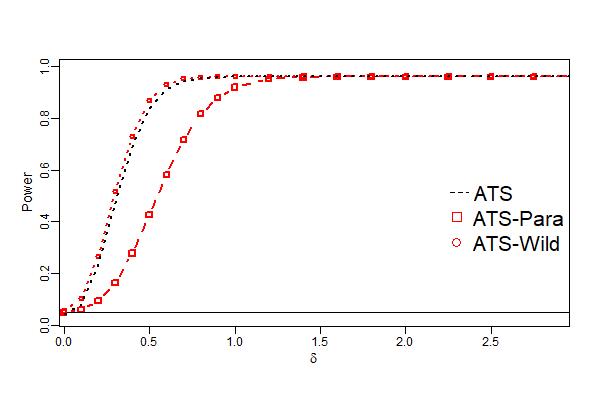}
\end{minipage}

\caption{{Simulated power for an one-point-alternative in scenario $B)$ ($\mathcal{H}_0^{\vv}:\vV_{1 11}=\vV_{2 11}=...=\vV_{1 55}$) { for the ATS tests based upon wild bootstrap, parametric bootstrap and Monte-Carlo critical values} as well as the two bootstrap tests based on Bartlett's statistic. The d=5 dimensional error terms are based on the skewed normal distribution with  covariance matrix  $(\vV)_{ij}=0.6^{|i-j|}$ and sample sizes $n_1=30 , n_2=20$ in the first row and $n_1=60, n_2=40$ in the second.}}
\label{Power3}

\end{figure}

\begin{figure}[htbp]
\begin{minipage}{\textwidth}
     \centering
\includegraphics[trim= 7mm 13mm 5mm 20mm,clip,scale=0.79]{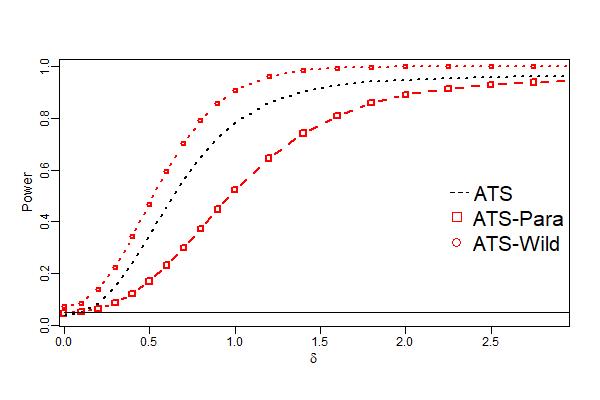}
 \end{minipage}
  \begin{minipage}{\textwidth}
     \centering
        \includegraphics[trim= 7mm 12mm 5mm 20mm,clip,scale=0.79]{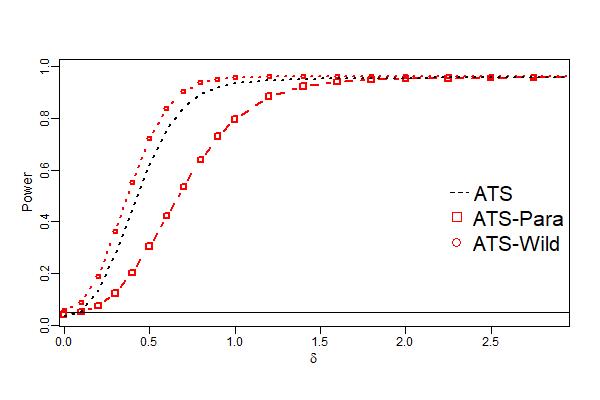}
\end{minipage}

\caption{{Simulated power for an one-point-alternative in scenario $B)$ ($\mathcal{H}_0^{\vv}:\vV_{1 11}=\vV_{2 11}=...=\vV_{1 55}$))  { for the ATS tests based upon wild bootstrap, parametric bootstrap and Monte-Carlo critical values} as well as the two bootstrap tests based on Bartlett's statistic. The d=5 dimensional error terms are based on the gamma  distribution with  covariance matrix  $(\vV)_{ij}=0.6^{|i-j|}$ and sample sizes $n_1=30 , n_2=20$ in the first row and $n_1=60, n_2=40$ in the second.}}\label{Power4}
\end{figure}
\begin{figure}[htbp]
\begin{minipage}{\textwidth}
     \centering
\includegraphics[trim= 7mm 13mm 5mm 20mm,clip,scale=0.79]{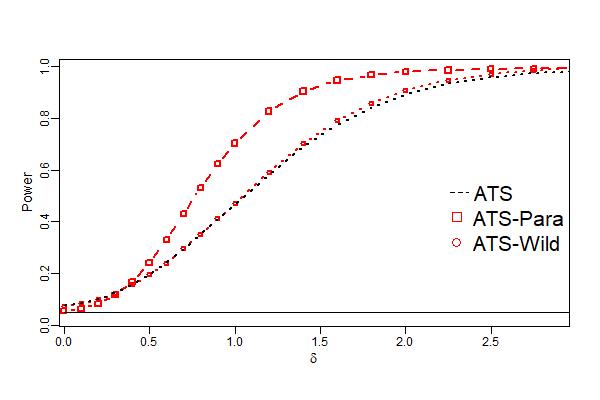}
 \end{minipage}
  \begin{minipage}{\textwidth}
     \centering
      \includegraphics[trim= 7mm 12mm 5mm 20mm,clip,scale=0.79]{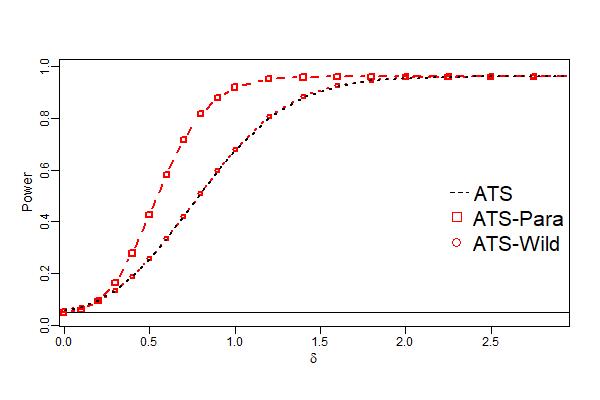}
\end{minipage}

\caption{{Simulated power for an one-point-alternative in scenario $C)$ ($\mathcal{H}_0^{\vv}:\tr(\vV_{1})=\tr(\vV_{2})$) { for the ATS tests based upon wild bootstrap, parametric bootstrap and Monte-Carlo critical values} as well as the two bootstrap tests based on Bartlett's statistic. The d=5 dimensional error terms are based on the skewed normal distribution with  covariance matrix  $(\vV)_{ij}=0.6^{|i-j|}$ and sample sizes $n_1=30 , n_2=20$ in the first row and $n_1=60, n_2=40$ in the second.}}
\label{Power5}

\end{figure}

\begin{figure}[H]
\begin{minipage}{\textwidth}
     \centering
\includegraphics[trim= 7mm 13mm 5mm 20mm,clip,scale=0.79]{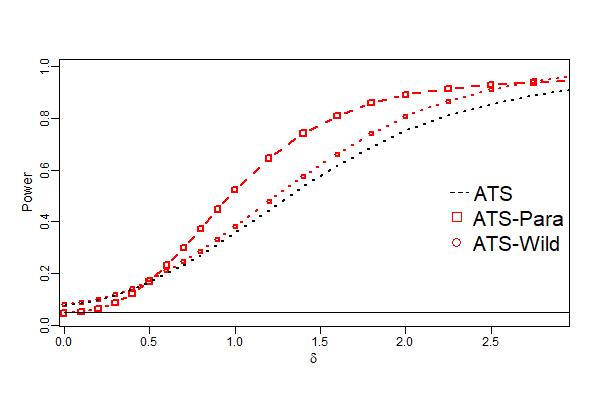}
 \end{minipage}
  \begin{minipage}{\textwidth}
     \centering
         \includegraphics[trim= 7mm 12mm 5mm 20mm,clip,scale=0.79]{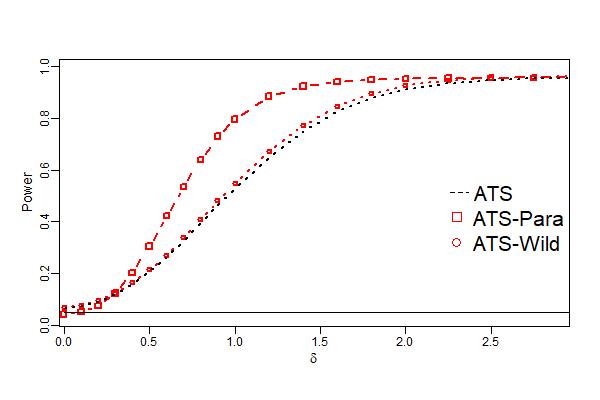}
\end{minipage}

\caption{{Simulated power for an one-point-alternative in scenario $C)$ ($\mathcal{H}_0^{\vv}:\tr(\vV_{1})=\tr(\vV_{2})$) { for the ATS tests based upon wild bootstrap, parametric bootstrap and Monte-Carlo critical values} as well as the two bootstrap tests based on Bartlett's statistic. The d=5 dimensional error terms are based on the gamma distribution with  covariance matrix  $(\vV)_{ij}=0.6^{|i-j|}$ and sample sizes $n_1=30 , n_2=20$ in the first row and $n_1=60, n_2=40$ in the second.}}
\label{Power6}

\end{figure}
\renewcommand{\baselinestretch}{2}

\subsection{Additional time computations }

 Here we present the time computation results for more hypotheses and test statistics as well as the required time in minutes for quadratic and for non-quadratic hypothesis matrices instead of
the relation between both. The considered hypotheses are:
\begin{itemize}
\item[$A$)] {Equal Covariance Matrices:\ } Testing $\mathcal{H}_0^{\vv}:\{\vV_{1}=\vV_{2}\} = \{ \vC(A) \vv = \vnull \}$ is usually described by $\vC(A)=\vP_2\otimes \vI_p$. However, the alternative choice $\widetilde\vC(A)=(1,-1)\otimes\vI_p\in \R^{p\times 2p}$ is computationally more efficient.

\item[$B$)] Equal Diagonal Elements: \ The hypothesis  $\mathcal{H}_0^{\vv}:\{\vV_{1 11}=...=\vV_{1 dd}\} = \{ \vC(B) \vv = \vnull \}$ can, e.g., be described by $\vC(B)=\diag(\vh_d)-\vh_d\cdot \vh_d^\top/d$. In contrast, the equivalent description by $\widetilde \vC(B)=(\veins_{d-1},\vnull_{(d-1)\times (d-1)},-\ve_1,\vnull_{(d-1)\times (d-2)},$ $-\ve_2,...,\vnull_{d-1},\ve_{d-1})\in \R^{(d-1)\times p}$ saves a considerable amount of time. Here, $\ve_j$ denotes the $d-1$ dimensional vector containing $1$ in the j-th component and $0$ elsewhere.
\item[$C$)] Equal traces: \ Testing the hypothesis  $\mathcal{H}_0^{\vv}:\{\tr(\vV_{1})=\tr(\vV_{2})\} = \{
\vC(C) \vv = \vnull\}$  is usually described by $\vC(C)= \vP_2\otimes [\vh_d\cdot \vh_d^\top]/d$. An equivalent expression is achieved with the smaller matrix $\widetilde\vC(C)=(1,-1)\otimes \vh_d/d\in \R^{1\times 2p}$.
\item[$D$)] Test for a given trace: \ $\mathcal{H}_0^{\vv}:\{\tr(\vV_{1})=\gamma\}\{
\vC(D) \vv = \vh_d \} $ for a given value $\gamma \in \R$ can either be described by  $\widetilde \vC(D)=\vh_d^\top/d\in \R^{1\times p}$ or $\vC(D)= [\vh_d\cdot \vh_d^\top]/d$, where the first choice has considerably less  rows. \end{itemize}

We used 4 different distributions (based on $t_9$-distribution, Normal-distribution, Skew Normal-distribution and Gamma-distribution)  and 2 covariance matrices ($(\vV_1)_{i,j}=0.6^{|i-j|}$ and $\vV_2=\vI_d+\vJ_d$) for each hypothesis and test-statistic. The average time of 100 such simulation runs are compared to get more valid results.\\
 For each test 1.000 bootstrap runs were performed with $n_1=125$ observations resp. $\vn=(150,100)$ observations in various dimensions. For the Monte-Carlo-test again 10.000 simulation steps are used.\\
The results confirm the remarks from \Cref{Time}. As expected, the WTS based on critical values on $\chi^2$ is even faster than the Monte-Carlo- ATS because no repetitions have to be done.
Moreover, it can be seen that both hypotheses for one group resp. two groups have comparable time demand for the quadratic idempotent hypothesis matrices. Therefore, here the concrete hypothesis seems to have no essential influence.
For the more time-efficient non-quadratic matrices, this is different, as it can be seen especially in comparison of $B)$ and $D)$. While the time for hypothesis $D)$ barely increases with the dimension,  there is a clear growth for $B)$. This fits with the number of lines of the according to hypothesis matrices.

\renewcommand{\baselinestretch}{1}
\begin{table}[htbp]

\centering

       \begin{tabular}{|l||r|r|r|r|l|r|r|r|r|}\hline

&\multicolumn{4}{|c|}{$C(A)$}&& \multicolumn{4}{|c|}{$\widetilde C(A)$}\\\cline{1-5}\cline{7-10}

  d   & 2&5&10&20  &      & 2&5&10&20   

       \\\cline{1-5}\cline{7-10}

ATS-Para & 0.757 & 5.401 & 27.928 & 222.443 && 0.745 & 5.332 & 27.227 & 195.237  \\\cline{1-5}\cline{7-10}

  ATS-Wild &0.451 & 0.612 & 10.831 & 114.175 && 0.455 & 0.599 & 10.034 & 87.011 \\ \cline{1-5}\cline{7-10}

  ATS &0.086 & 0.195 & 0.387 & 1.231 && 0.056 & 0.154 & 0.276 & 0.698    \\ \cline{1-5}\cline{7-10}

  WTS-Para &0.869 & 6.625 & 44.295 & 462.852 && 0.839 & 6.157 & 34.851 & 275.909 \\ \cline{1-5}\cline{7-10}

  WTS-Wild  &0.559 & 0.916 & 27.252 & 355.679 && 0.545 & 0.794 & 17.681 & 170.063 \\ \cline{1-5}\cline{7-10}

   WTS-$\chi^2$ & 0.003 & 0.004 & 0.038 & 0.320 && 0.003 & 0.003 & 0.029 & 0.133 \\

   \hline

\end{tabular}

\caption{ Required time in seconds for various tests statistics and different dimensions for hypothesis $A)$ ($\mathcal{H}_0^{\vv}:\vV_{1}=\vV_{2}$) with a quadratic hypothesis matrix on the left side and a non-quadratic hypothesis matrix on the right sight.}

\end{table}

\begin{table}[htbp]

\centering

       \begin{tabular}{|l||r|r|r|r|l|r|r|r|r|}\hline

&\multicolumn{4}{|c|}{$C(B)$}&& \multicolumn{4}{|c|}{$\widetilde C(B)$}\\\cline{1-5}\cline{7-10}

  d   & 2&5&10&20  &      & 2&5&10&20

       \\\cline{1-5}\cline{7-10}

ATS-Para & 0.441 & 2.795 & 11.367 & 76.967 && 0.374 & 0.454 & 1.152 & 3.380   \\ \cline{1-5}\cline{7-10}

 ATS-Wild &0.285 & 0.349 & 0.743 & 5.159 && 0.275 & 0.287 & 0.316 & 0.379  \\ \cline{1-5}\cline{7-10}

 ATS& 0.048 & 0.152 & 0.258 & 0.594 && 0.027 & 0.058 & 0.111 & 0.162\\ \cline{1-5}\cline{7-10}

 WTS-Para& 0.534 & 3.390 & 15.386 & 120.942 && 0.454 & 0.551 & 1.411 & 4.499 \\ \cline{1-5}\cline{7-10}

 WTS-Wild& 0.364 & 0.491 & 8.341 & 74.507 && 0.342 & 0.374 & 0.444 & 0.661 \\ \cline{1-5}\cline{7-10}

  WTS-$\chi^2$& 0.002 & 0.002 & 0.018 & 0.090 && 0.001 & 0.002 & 0.003 & 0.003  \\

   \hline

   \hline

\end{tabular}

\caption{ Required time in seconds for various tests statistics and different dimensions for hypothesis $B)$ ($\mathcal{H}_0^{\vv}:\vV_{1 11}=...=\vV_{1 dd}$)  with a quadratic hypothesis matrix on the left side and a non-quadratic hypothesis matrix on the right sight.}

\end{table}

\begin{table}[htbp]

\centering

       \begin{tabular}{|l||r|r|r|r|l|r|r|r|r|}\hline

&\multicolumn{4}{|c|}{$C(C)$}&& \multicolumn{4}{|c|}{$\widetilde C(C)$}\\\cline{1-5}\cline{7-10}

  d   & 2&5&10&20  &      & 2&5&10&20   

       \\\cline{1-5}\cline{7-10}

ATS-Para &0.751 & 5.419 & 27.971 & 222.583 && 0.735 & 5.263 & 26.721 & 178.510  \\\cline{1-5}\cline{7-10}

  ATS-Wild &0.451 & 0.607 & 10.736 & 114.369 && 0.446 & 0.579 & 9.394 & 70.564\\ \cline{1-5}\cline{7-10}

  ATS & 0.087 & 0.194 & 0.386 & 1.211 && 0.037 & 0.040 & 0.043 & 0.117  \\ \cline{1-5}\cline{7-10}

  WTS-Para &0.856 & 6.029 & 39.253 & 416.780 && 0.820 & 5.581 & 27.059 & 179.230  \\ \cline{1-5}\cline{7-10}

  WTS-Wild &  0.548 & 0.769 & 21.907 & 305.647 && 0.525 & 0.649 & 9.861 & 71.164 \\ \cline{1-5}\cline{7-10}

   WTS-$\chi^2$ &0.003 & 0.003 & 0.035 & 0.269 && 0.003 & 0.004 & 0.021 & 0.034 \\

   \hline

\end{tabular}

\caption{ Required time in seconds for various tests statistics and different dimensions for hypothesis $C)$ ($\mathcal{H}_0^{\vv}:\tr(\vV_{1})=\tr(\vV_{2})$)  with a quadratic hypothesis matrix on the left side and a non-quadratic hypothesis matrix on the right sight.}

\end{table}

\begin{table}[htbp]

\centering

       \begin{tabular}{|l||r|r|r|r|l|r|r|r|r|}\hline

&\multicolumn{4}{|c|}{$C(D)$}&& \multicolumn{4}{|c|}{$\widetilde C(D)$}\\\cline{1-5}\cline{7-10}

  d   & 2&5&10&20  &      & 2&5&10&20

       \\\cline{1-5}\cline{7-10}

 ATS-Para\hspace*{-0.1cm} & 0.443 & 2.700 & 11.133 & 76.557 && 0.377 & 0.376 & 0.374 & 0.402  \\ \cline{1-5}\cline{7-10}

 ATS-Wild \hspace*{-0.1cm} &0.286 & 0.350 & 0.745 & 5.136 && 0.275 & 0.275 & 0.274 & 0.302 \\ \cline{1-5}\cline{7-10}

 ATS& 0.049 & 0.152 & 0.258 & 0.589 && 0.028 & 0.029 & 0.029 & 0.032  \\ \cline{1-5}\cline{7-10}

WTS-Para\hspace*{-0.1cm} &  0.534 & 3.111 & 14.416 & 118.989 && 0.452 & 0.455 & 0.456 & 0.476\\ \cline{1-5}\cline{7-10}

 WTS-Wild\hspace*{-0.1cm} & 0.363 & 0.452 & 6.870 & 73.208 && 0.342 & 0.344 & 0.341 & 0.370  \\ \cline{1-5}\cline{7-10}

 WTS-$\chi^2$& 0.002 & 0.002 & 0.016 & 0.089 && 0.001 & 0.002 & 0.002 & 0.021 \\

   \hline

\end{tabular}

\caption{ Required time in seconds for various tests statistics and different dimensions for hypothesis $D)$ $\mathcal{H}_0^{\vv}:\{\tr(\vV_{1})=\gamma\}$  with a quadratic hypothesis matrix on the left side and a non-quadratic hypothesis matrix on the right sight.}

\end{table}

\bibliographystyle{apalike}
\bibliographystyle{unsrtnat}
\bibliography{Literatur}
\vskip .65cm
\noindent
TU Dortmund University, Faculty of Statistics, Germany
\vskip 2pt
\noindent
E-mail: (email: paavo.sattler@tu-dortmund.de)
\vskip 2pt

\noindent
Department for Mathematics, University of Salzburg, Austria
\vskip 2pt
\noindent
E-mail: (email: Arne.Bathke@sbg.ac.at)
\vskip 2pt

\vskip 2pt
TU Dortmund University, Faculty of Statistics, Germany
\vskip 2pt
\noindent
E-mail: (email: markus.pauly@tu-dortmund.de)
\vskip 2pt
\end{document}